\documentclass[12pt,a4paper,reqno]{amsart}
\usepackage{color}
\usepackage{graphics,epic}
\usepackage{amsmath,amssymb, amsthm}
\usepackage[all,2cell]{xy}

\newtheorem{theorem}{Theorem}[section]
\newtheorem*{theorem*}{Theorem}

\newtheorem{lemma}[theorem]{Lemma}
\newtheorem{proposition}[theorem]{Proposition}

\newtheorem*{conjecture*}{Conjecture}

\newtheorem{remark}[theorem]{Remark}

\renewcommand{\hat}[1]{\widehat{#1}}

\newcommand{\id}{{\rm id}}

\newcommand{\im}{{\rm im}}

\newcommand{\End}{{\rm End}\,}
\newcommand{\Res}{{\rm Res}\,}

\newcommand{\Z}{\mathbb{Z}}
\newcommand{\Q}{\mathbb{Q}}
\newcommand{\C}{\mathbb{C}}

\def\Res{{\rm Res}}
\def\wt{{\rm wt}}
\def\C{{\mathbb C}}

\def\R{{\mathbb R}}

\def\Z{{\mathbb Z}}

\def\1{{\bf 1}}

\def \End{{\rm End}}

\def \Ind{{\rm Ind}}

\def \pf{\noindent {\bf Proof: \,}}
\def\theequation{5.\arabic{equation}}
\def \h{\mathfrak{h}}
\def \w{\omega}
\def \g{\mathfrak{g}}

\begin{document}

\title[Rationality of VOSAs with rational conformal weights]{Rationality of vertex operator superalgebras with rational conformal weights}


\author{Xingjun Lin}
\address{Xingjun Lin,  School of Mathematics and Statistics, Wuhan University, Wuhan 430072,  China.}
\thanks{X. Lin was supported by China NSF grants
11801419, 12171371 and the starting research fund from Wuhan University}
\email{linxingjun88@126.com}
\begin{abstract}
For the affine vertex algebra $V_k(\g)$ at an admissible level $k$ of $\hat{\g}$, we prove that certain subcategory of weak  $V_k(\g)$-module category is semisimple. As a consequence, we show that  $V_k(\g)$ is rational with respect to a family of Virasoro elements. We also prove that certain affine vertex operator superalgebras and minimal $W$-algebras are rational with respect to a family of Virasoro elements. 
\end{abstract}
\maketitle
\section{Introduction}\label{intro}
\def\theequation{1.\arabic{equation}}
\setcounter{equation}{0}

Rationality and $C_2$-cofiniteness are important properties of  vertex operator algebras. Rational and $C_2$-cofinite vertex operator algebras have several good properties. For instance, for a rational and $C_2$-cofinite vertex operator algebra $V$, the space of trace functions of irreducible $V$-modules is invariant under the action of the modular group  $SL(2, \Z)$ \cite{Z}, \cite{DLM2}. Another good property of  rational and $C_2$-cofinite vertex operator algebras is that the categories of their modules are modular tensor categories \cite{Hu}.

Affine vertex operator algebras at positive integer levels provide important examples of rational and $C_2$-cofinite vertex operator algebras \cite{FZ}, \cite{DLM2}. Moreover, it was proved in  \cite{DM} that an affine vertex operator algebra is rational and $C_2$-cofinite if and only if its level is a positive integer.  On the other hand, Kac-Wakimoto proved in \cite{KW}, \cite{KW1} that the space of characters of irreducible highest weight modules of an affine Lie algebra at a fixed admissible level is invariant under the action of the modular group  $SL(2, \Z)$.  However, affine vertex operator algebras at admissible levels may not be rational \cite{AM}, \cite{DLM0}.

It is a natural problem to give an explanation of the results of  Kac-Wakimoto in the framework of vertex operator algebras. To solve this problem, Van Ekeren established modular invariance properties of rational $\Q$-graded vertex operator superalgebras in \cite{E}. Thus, it is important to prove that affine vertex algebras at admissible levels are rational $\Q$-graded vertex operator algebras. In the special case that the Lie algebra is $sl_2$, this has been established in \cite{DLM0}. In the general case, it was proved in \cite{AM}, \cite{A1} that affine vertex algebras at admissible levels are rational in the BGG category. Based on these results, Arakawa-Van Ekeren explained  the results of  Kac-Wakimoto in the framework of vertex operator algebras in  \cite{AE}.

For a $\Q$-graded vertex operator algebra $V$, we may define weak modules, admissible modules and ordinary modules of $V$ (cf. Subsection \ref{bas}). It was proved in \cite{AE} that the categories of ordinary modules of affine vertex operator algebras at admissible levels are semisimple. To establish the modular invariance property of a $\Q$-graded vertex operator superalgebra $V$,  one of key conditions is that the Zhu's algebra of $V$ is semisimple. For a $\Q$-graded vertex operator algebra $V$, it is not known whether the Zhu's algebra  of $V$ is finite dimensional. Thus, to obtain naturally semisimplicity of the Zhu's algebra, we define a $\Q$-graded vertex operator algebra $V$ to be rational if the category of admissible modules of $V$ is semisimple.  In this paper, we show that affine vertex operator algebras at admissible levels are rational with respect to a family of Virasoro elements (see Theorem \ref{rationalaff}).

It is a natural problem to find more examples of rational $\Q$-graded vertex operator superalgebras. In this paper, we show that certain vertex superalgebras associated with affine Lie superalgebras are rational $\Q$-graded vertex operator superalgebras (see Theorems  \ref{r-B(0,n)}, \ref{r-s-2}). We also prove that certain minimal $W$-algebras are rational $\Q$-graded vertex operator algebras (see Theorems \ref{r-w-1}, \ref{r-w-2}). The reason to consider these vertex operator superalgebras is that they are extensions of affine vertex operator algebras at admissible levels \cite{AMPP}, \cite{AKMPP1}, \cite{AKMPP2}. To prove rationality of these vertex operator superalgebras, we show that certain categories  of weak modules of  extensions of vertex operator algebras are semisimple (see Theorem \ref{scextension}).

The paper is organized as follows: In Section 2, we recall basic definitions about vertex operator superalgebras and basic facts about rational vertex operator superalgebras. In Section 3, we prove that affine vertex operator algebras at admissible levels are rational with respect to a family of Virasoro elements. One of our goals is to show that certain affine vertex operator superalgebras and minimal $W$-algebras are rational with respect to a family of  Virasoro elements. To obtain these results,   we prove that certain categories  of weak modules of  extensions of vertex operator algebras are semisimple in Section 4. In Section \ref{affineVOSA}, we prove that certain affine vertex operator superalgebras are rational with respect to a family of Virasoro elements. In Section 6, we show that certain minimal $W$-algebras are rational with respect to a family of Virasoro elements. 

\section{Preliminaries }
\def\theequation{2.\arabic{equation}}
\setcounter{equation}{0}
\subsection{Basics}\label{bas}
In this subsection, we recall from \cite{K2}, \cite{KWang}, \cite{L1}, \cite{E} some facts about vertex operator superalgebras with rational conformal weights. Let $V=V_{\bar 0}\oplus V_{\bar 1}$ be a vector superspace, the element in $V_{\bar 0} $ (resp. $V_{\bar 1}$) is called  {\em even} (resp. {\em odd}).
For any $v\in V_{\bar i}$ with  $i=0,1$, we define the parity of $v$ to be $[v]=i$.  A {\em vertex superalgebra} is a quadruple $(V,Y(\cdot, z),\1,D),$ where $V=V_{\bar 0}\oplus V_{\bar 1}$ is a vector superspace, ${\bf 1}$ is an even vector of $V$ called the {\em vacuum vector} of $V$, $D$ is an endomorphism of $V$,   and $Y(\cdot, z)$ is a parity preserving linear map
\begin{align*}
 Y(\cdot, z): V &\to (\End\,V)[[z,z^{-1}]] ,\\
 v&\mapsto Y(v,z)=\sum_{n\in{\Z}}v_{(n)}z^{-n-1}\ \ \ \  (v_{(n)}\in
\End\,V)
\end{align*}
satisfying the following axioms:  \\
(i) For any $u,v\in V,$ $u_{(n)}v=0$ for sufficiently large $n$;\\
(ii) $Y({\bf 1},z)=\id_{V}$;\\
(iii) $Y(v,z){\bf 1}=v+\sum_{n\geq 2}v_{(-n)}{\bf 1}z^{n-1},$ for any $v\in V$;\\
(iv) $[D, Y(v, z)]=Y(D(v), z)=\frac{d}{dz}Y(v, z)$;\\
(v) The {\em Jacobi identity} for ${\Z}_2$-homogeneous $u,v\in V$ holds,
\begin{align*}
\begin{array}{c}
\displaystyle{z^{-1}_0\delta\left(\frac{z_1-z_2}{z_0}\right)
Y(u,z_1)Y(v,z_2)-(-1)^{[u][v]} z^{-1}_0\delta\left(\frac{z_2-z_1}{-z_0}\right)
Y(v,z_2)Y(u,z_1)}\\
\displaystyle{=z_2^{-1}\delta
\left(\frac{z_1-z_0}{z_2}\right)
Y(Y(u,z_0)v,z_2)}.
\end{array}
\end{align*}
This completes the definition of a vertex superalgebra and we will denote the vertex superalgebra briefly by $V$.

A vertex superalgebra $V$ is called a {\em vertex operator superalgebra} if there is an even vector  $\omega$ called the {\em Virasoro element} of $V$ such that the following two conditions hold: \\
 (vi) The component operators of  $Y(\w,z)=\sum_{n\in\Z}L(n)z^{-n-2}$ satisfy the Virasoro algebra
relation with {\em central charge} $c\in \C:$
\begin{align*}
[L(m),L(n)]=(m-n)L(m+n)+\frac{1}{12}(m^3-m)\delta_{m+n,0}c,
\end{align*}
and
$$L(-1)=D;$$
(vii) $V$ is $\Q$-graded such that $V=\oplus_{n\in \Q}V_n$,  $L(0)|_{V_n}=n$, $\dim(V_n)<\infty$ and $V_n=0$ for sufficiently small $n$. For $v\in V_n$, the {\em conformal weight} of $v$ is defined to be $\wt u=n$.

\begin{remark}
In this paper, the conformal weights of a vertex operator superalgebra could be rational numbers. Usually, the conformal weights of a vertex operator superalgebras are required to be integers or half-integers (cf. \cite{DZ}, \cite{DZ1}, \cite{FLM}, \cite{LL}, \cite{KWang}). In this paper, vertex operator algebras having integral conformal weights are referred to as $\Z$-graded vertex operator algebras.
\end{remark}

For a vertex operator superalgebra, we may define three classes of modules. First, let $V$ be a vertex superalgebra. A {\em weak  $V$-module} is a vector superspace $M=M_{\bar 0}\oplus M_{\bar 1}$ equipped
with a parity preserving linear map
\begin{align*}
Y_{M}(\cdot, z):V&\to (\End M)[[z, z^{-1}]],\\
v&\mapsto Y_{M}(v,z)=\sum_{n\in\Z}v_{(n)}z^{-n-1},\,v_{(n)}\in \End M
\end{align*}
satisfying the following conditions: For any ${\Z}_2$-homogeneous $u,v\in V$, $w\in M$ and $n\in \Z$,
\begin{align*}
&\ \ \ \ \ \ \ \ \ \ \ \ \ \ \ \ \ \ \ \ \ \ \ \ \ \ u_{(n)}w=0 \text{ for  sufficiently large } n;\\
&\ \ \ \ \ \ \ \ \ \ \ \ \ \ \ \ \ \ \ \ \ \ \ \ \ \ Y_M(\1, z)=\id_M;\\
&z^{-1}_0\delta\left(\frac{z_1-z_2}{z_0}\right)
Y_M(u,z_1)Y_M(v,z_2)w-(-1)^{[u][v]} z^{-1}_0\delta\left(\frac{z_2-z_1}{-z_0}\right)
Y_M(v,z_2)Y_M(u,z_1)w\\
&\ \ \ \ \ \ \ \ \ \ \ \ \ \ \ \ =z_2^{-1}\delta
\left(\frac{z_1-z_0}{z_2}\right)
Y_M(Y(u,z_0)v,z_2)w.
\end{align*}


\vskip.5cm
We next define the second class of modules. Let $V$ be a vertex superalgebra, and $\omega$ be a Virasoro element of $V$.
A weak
 $V$-module $M$  is called an \textit{$\omega$-admissible $V$-module} if $M$  has an $\R$-gradation $M=\bigoplus_{n\in \R, n\geq 0}M(n)$ such
that
\begin{align*}\label{AD1}
a_{(m)}M(n)\subset M(\wt a+n-m-1)
\end{align*}
for any homogeneous $a\in V$, $m\in\Z$ and $n\in \R$.
An $\omega$-admissible $V$-module $M$ is said to be
\textit{irreducible} if $M$ has no non-trivial graded weak
$V$-submodule.
A vertex superalgebra $V$ is called \textit{$\omega$-rational} if
any  $\omega$-admissible $V$-module is a direct sum of irreducible $\omega$-admissible $V$-modules.
\begin{remark}
For a vertex superalgebra $V$, there may be many Virasoro elements such that $V$ is a vertex operator superalgebra. We will see many examples of  vertex superalgebras which are rational with respect to some Virasoro elements, but there  exist also Virasoro elements such that the vertex superalgebras are not rational.
\end{remark}

We next define the third class of modules. Let $V$ be a vertex superalgebra, and $\omega$ be a Virasoro element of $V$. A weak $V$-module $M$ is called an {\em $\omega$-ordinary $V$-module} if  $M=\bigoplus_{\lambda\in\C}
M_{\lambda}$ such that $M_\lambda=\{w\in M|L(0)w=\lambda w\}$, $M_\lambda$ is
finite dimensional and $M_{\lambda}=0$ if the real part of $\lambda$ is sufficiently small.

The relation between the three classes of modules of $\Z$-graded vertex operator algebra has been discussed in \cite{DLM3}. By the similar argument as in the proof of Lemma 3.4 of \cite{DLM3}, we have the following result.
\begin{proposition}\label{ordi}
Let $V$ be a vertex superalgebra, $\omega$ be a Virasoro element of $V$ such that the eigenvalues of $\omega_{(1)}$ are contained in $\frac{1}{T}\Z$ for some positive integer $T$. If  $M$ is an $\omega$-ordinary $V$-module, then $M$ is an $\omega$-admissible $V$-module.
\end{proposition}

\begin{remark}
 (1) Our definition of $\omega$-rationality is different from that in Definition 2.6 of \cite{AE}. By Definition 2.6 of \cite{AE}, a vertex operator superalgebra $(V,Y(\cdot, z),\1,\omega)$ is rational if $V$ has finitely many irreducible $\omega$-ordinary $V$-modules and every $\omega$-ordinary $V$-module is a direct sum of irreducible $\omega$-ordinary $V$-modules. Under the assumption of Proposition \ref{ordi}, if $V$ is $\omega$-rational then $V$ is rational in the sense of Arakawa-Van Ekeren \cite{AE}.

 (2)  For a $\Q$-graded vertex operator algebra $V$, it is not known whether the Zhu's algebra  of $V$ is finite dimensional. Thus, the assumption that  every $\omega$-ordinary $V$-module is a direct sum of irreducible $\omega$-ordinary $V$-modules does not imply that the Zhu's algebra of $V$  is semisimple. To obtain naturally semisimplicity of the Zhu's algebra, it is necessary to consider $\omega$-rational vertex operator superalgebras.
  \end{remark}

We now recall another important notion from \cite{AE}. Let $V^0$ be a vertex superalgebra, $V=V^0\oplus V^+$ be an extension of $V^0$ by its weak module $V^+$. Set
\begin{align*}
&C^{rel}(V)=V^0_{(-2)}V^0+V_{(-1)}V^+=V_{(-2)}V+V_{(-1)}V^+,\\
&R^{rel}(V)=V/C^{rel}(V),
  \end{align*}where $U_{(-1)}W=\{u_{(-1)}w|u\in U, w\in W\}$ and $U_{(-2)}W=\{u_{(-2)}w|u\in U, w\in W\}$ for subspaces $U, W$ of $V$. Then we say $V$ is {\em cofinite} relative to the decomposition $V=V^0\oplus V^+$ if $\dim R^{rel}(V)<\infty$.
  \begin{remark}
  In Proposition 5.12 and Theorem 5.13 of \cite{AE}, Arakawa and Van Ekeren showed that vertex superalgebras have modular invariance properties under certain assumptions about $\omega$-rationality and relative cofiniteness. Hence,  it is important to study $\omega$-rationality and relative cofiniteness of vertex superalgebras.
  \end{remark}
  \subsection{Change of Virasoro elements} In this subsection, we recall from \cite{AE}, \cite{DLinM} a way to change the Virasoro element. Let $(V,Y(\cdot, z),\1,\omega)$ be a $\Z$-graded vertex operator algebra, $h\in V_1$ be a vector satisfying the following conditions:
  $$[h_{(m)}, h_{(n)}]=2m\delta_{m+n, 0},$$
  $$[L(m), h_{(n)}]=-nh_{(m+n)}-\frac{m^2+m}{2}\delta_{m+n, 0}p.$$
  Furthermore, we assume that $h_{(0)}$ acts semisimply on $V$ and  the eigenvalues of $h_{(0)}$ are rational numbers. Since $L(0)$ and $h_{(0)}$ are commutative, each $V_n$ is a direct sum of eigenspaces of $h_{(0)}$. For a rational number $\sigma$, set $\omega(\sigma)=\omega-\frac{\sigma}{2}L(-1)h$. It is known that the component operators of $\omega(\sigma)$ satisfy the Virasoro algebra relations (cf. \cite{AE}). We will see that  $(V,Y(\cdot, z),\1,\omega(\sigma))$ is a vertex operator algebra under certain assumption. Set $$V_{m, n}=\{v\in V_m|h_{(0)}v=nv\}, $$ for $m\in \Z$ and $n\in \Q$. Then $$V'_{m}=\{v\in V|(\omega-\frac{\sigma}{2}L(-1)h)_{(1)}v=mv\}=\coprod_{s\in \Z, t\in \Q, s+\frac{\sigma t}{2}=m}V_{s, t}.$$  The following result was essentially proved in  \cite{DLinM}, \cite{AE}.
  \begin{proposition}\label{Virasoro}
  Let $(V,Y(\cdot, z),\1,\omega)$ be a $\Z$-graded vertex operator algebra of central charge $c$, $h\in V_1$ be a vector satisfying the conditions above. Suppose that  $ V_{m}'$ is finite dimensional for any $m\in \Q$ and  $ V_{m}'=0$ if $m$ is sufficiently small. Then   $(V,Y(\cdot, z),\1,\omega(\sigma))$ is a vertex operator algebra of central charge $c+6\sigma(p-\sigma)$.
  \end{proposition}

\subsection{Regular vertex operator algebras}
In this subsection, we recall some facts about regular vertex operator algebras from \cite{DLM1} and \cite{DY}.
 A $\Z$-graded vertex operator algebra $(V,Y(\cdot, z),\1,\omega)$ is said to be {\em regular} if any weak $V$-module $M$ is a direct sum of irreducible $\omega$-ordinary $V$-modules. The following result has been proved in Corollary 4.6 of \cite{DY}.
 \begin{proposition}\label{regu}
 Let $(V,Y(\cdot, z),\1,\omega)$ be a $\Z$-graded vertex operator algebra. Then the following statements are equivalent:\\
 (1) The category of weak $V$-modules is semisimple.\\
 (2) $V$ is regular.
 \end{proposition}
As a consequence, assume that $(V,Y(\cdot, z),\1,\omega)$ is a regular $\Z$-graded vertex operator algebra and $\omega'$ is another Virasoro element of $V$. By Proposition \ref{regu}, $V$ is $\omega'$-rational.
 \subsection{Weyl vertex algebras}
 In this subsection, we will show that the Weyl vertex algebra is rational with respect to a family of Virasoro elements.

 First, we recall from \cite{K2}, \cite{AP} some facts about the Weyl vertex algebra. The Weyl algebra $\hat{\mathcal A}$ is an associative algebra with generators
 $$a_n, a_n^*~(n\in \Z),$$
 and relations
 $$[a_n, a^*_m]=\delta_{m+n, 0}, [a_n, a_m]=0=[a^*_m, a^*_n] ~(n, m\in \Z).$$
 Let $W$ denote the simple Weyl module generated by the cyclic vector $\1$ such that $$a_n\1=a^*_{n+1}\1=0~(n\geq 0).$$
 As a vector space,
 $$W\cong \C[a_{-n}, a^*_{-m}|n>0, m\geq 0].$$
 There is a unique vertex algebra $(W, Y(\cdot, z), \1)$ (cf. \cite{K2}, \cite{AP}) such that
 $$Y(a_{-1}\1, z)=\sum_{n\in \Z}a_nz^{-n-1},~Y(a^*_{0}\1, z)=\sum_{n\in \Z}a^*_nz^{-n}.$$
 Set $\omega=a_{-1}a^*_{-1}\1$, then it was known that the component operators of $\omega$ satisfy  the  Virasoro relations (cf. \cite{K2}, \cite{AP}).

 Set $\beta=a_{-1}a^*_{0}\1$, then it was known \cite{K2}, \cite{AP} that $[\beta_{(n)}, \beta_{(m)}]=-n\delta_{n+m, 0}$, and $[\beta_{(n)}, a_m]=-a_{n+m},~[\beta_{(n)}, a^*_m]=a^*_{n+m}$.
 As a consequence, we have the following result which was essentially obtained in \cite{K2}, \cite{AP}.
  \begin{proposition}
  Let $\mu$ be a rational number such that $0<\mu<1$, and set $\omega_{\mu}=\omega-\mu\beta_{(-2)}\1$. Then $(W, Y(\cdot, z), \1, \omega_{\mu})$ is a vertex operator algebra. Moreover,  the conformal weights of $a_{-1}\1$, $a^*_{0}\1$ are $1-\mu$, $\mu$, respectively.
   \end{proposition}
   \pf By Proposition \ref{Virasoro}, the component operators of $\omega_{\mu}$ satisfy  the  Virasoro relations. Moreover, it was known  \cite{AP} that $$[\omega_{(1)}, a_m]=-ma_{m},~[\omega_{(1)}, a^*_{m}]=-ma^*_m.$$
   Hence, $$[(\omega_{\mu})_{(1)}, a_m]=(-\mu-m)a_{m},~[(\omega_{\mu})_{(1)}, a^*_{m}]=(\mu-m)a^*_m.$$
   Since $W\cong \C[a_{-n}, a^*_{-m}|n>0, m\geq 0],$   $(W, Y(\cdot, z), \1, \omega_{\mu})$ is a vertex operator algebra by Proposition \ref{Virasoro}.
   \qed
   \begin{proposition}
  Let $\mu$ be a rational number such that $0<\mu<1$. Then $(W, Y(\cdot, z), \1)$ is $ \omega_{\mu}$-rational.
   \end{proposition}
   \pf Define the operators $a_{\{n\}}, a^*_{\{n\}} ~(n\in \frac{1}{2}+\Z)$ by
   $$Y(a_{-1}\1, z)=\sum_{n\in \frac{1}{2}+\Z}a_{\{n\}}z^{-n-1/2},~Y(a^*_{0}\1, z)=\sum_{n\in \frac{1}{2}+ \Z}a^*_{\{n\}}z^{-n-1/2}.$$
   Then we have $$[a_{\{n\}}, a^*_{\{m\}}]=\delta_{m+n, 0}, [a_{\{n\}}, a_{\{m\}}]=0=[a^*_{\{n\}}, a^*_{\{m\}}] ~(n, m\in \frac{1}{2}+\Z).$$
   Thus the operators $a_{\{n\}}, a^*_{\{n\}} ~(n\in \frac{1}{2}+\Z)$ satisfy the relation (3.6.1) of \cite{K2}. Set $\hat{\mathcal A}^+=\{a_{\{n\}}, a^*_{\{n\}}|n\geq \frac{1}{2}\}$. Let $M$ be a weak $W$-module such that $\hat{\mathcal A}^+$ acts locally nilpotently on $M$, i.e., for any  $m\in M$ there exists $n>0$ such that $g_1\cdots g_nm=0$ for any $n$ elements $g_1, \cdots, g_n$ of $\hat{\mathcal A}^+$. Then, by Theorem 3.6 of \cite{K2}, $M$ is a direct sum of copies of $W$.

   Note that $a^*_{\{n\}}=(a^*_{0}\1)_{(n-\frac{1}{2})}$,  $a_{\{n\}}=(a_{-1}\1)_{(n-\frac{1}{2})}$. As a consequence,   $\hat{\mathcal A}^+$ acts locally nilpotently on all the $\omega_{\mu}$-admissible $W$-modules. This implies $(W, Y(\cdot, z), \1)$ is $ \omega_{\mu}$-rational.
   \qed
\subsection{Zhu's algebras}
In this subsection, we will show that if  the vertex operator superalgebra $(V,Y(\cdot, z),\1,\omega)$ is $\omega$-rational, then the Zhu's algebra of $(V,Y(\cdot, z),\1,\omega)$ is semisimple.

First, we recall from \cite{DK}, \cite{DLM0}, \cite{Z} some facts about the Zhu's algebras of vertex operator superalgebras. Let $(V,Y(\cdot, z),\1,\omega)$ be a vertex operator superalgebra. For any $v\in V$, let $\gamma_v$ be the largest integer less than or equal to $\wt v$. Set $\epsilon_v=\gamma_v-\wt v$, then we define $$\chi(u,v)= \left\{
\begin{array}{ll}
1, & \text{ if } \epsilon_u+\epsilon_v\leq -1,\\
0,&  \text{ otherwise },
\end{array}\right.$$
for any $u, v\in V$.

For $u, v\in V$ such that $\epsilon_u+\epsilon_v\in \Z$, we define $u\circ_{-2+\chi(u, v)}v$ as follows:
$$u\circ_{-2+\chi(u, v)}v= \left\{
\begin{array}{ll}
\Res_z\frac{(1+z)^{\gamma_u}}{z^2}Y(u, z)v, & \text{ if } \epsilon_u+\epsilon_v=0,\\
\Res_z\frac{(1+z)^{\gamma_u}}{z}Y(u, z)v,&  \text{ if } \epsilon_u+\epsilon_v=-1,
\end{array}\right.$$
We then define $J$ to be the $\C$-span of the following set of elements
$$\{u\circ_{-2+\chi(u, v)}v|u, v\in V \text{ and } \epsilon_u+\epsilon_v\in \Z\}.$$
Set $V_{\Z}=\{v\in V|\wt v\in \Z\}$ and $Zhu_{\omega}(V)=V_{\Z}/J$. Then it is known \cite{DK} that $Zhu_{\omega}(V)$ is an associative superalgebra with the product defined as follows: For any $\bar{u}, \bar{v}\in Zhu_{\omega}(V)$,
$$\bar{u}\cdot\bar{v}=\overline{\Res_{z}\frac{(1+z)^{\gamma_u}}{z}Y(u, z)v}.$$
Moreover, the following results were established in \cite{Z}, \cite{KWang}, \cite{DK}.
\begin{theorem}\label{Zhu}
(1) There is a restriction functor $\Omega_{\omega}$ from the category of $\omega$-admissible V-modules to the category of $Zhu_{\omega}(V)$-modules. It sends $M$ to $M(0)$ with the action $\bar u\cdot x=u_{(\wt u-1)}x$ for $\bar{u}\in Zhu_{\omega}(V)$ and $x\in M(0)$.\\
(2) There is an induction functor $L_{\omega}$ from the category of $Zhu_{\omega}(V)$-modules to the category of $\omega$-admissible V-modules. Moreover, $\Omega_{\omega}(L_{\omega}(N))=N$ for any $Zhu_{\omega}(V)$-module $N$.\\
(3) $\Omega_{\omega}$ and $L_{\omega}$ are inverse bijections between the sets of irreducible modules in each category.
\end{theorem}

As a consequence, we have the following result which is the generalization of Theorem 8.1 of \cite{DLM3}.
\begin{proposition}
Let $(V,Y(\cdot, z),\1,\omega)$  be a vertex operator algebra. Suppose that $V$ is $\omega$-rational. Then the Zhu's algebra $Zhu_{\omega}(V)$ is  a finite dimensional  semisimple associative algebra.
\end{proposition}
\pf By Theorem \ref{Zhu}, $L_{\omega}(Zhu_{\omega}(V))$ is an $\omega$-admissible $V$-module. By assumption, $V$ is $\omega$-rational, this implies that $L_{\omega}(Zhu_{\omega}(V))$  is completely reducible. Since $\Omega_{\omega}(L_{\omega}(Zhu_{\omega}(V)))=Zhu_{\omega}(V)$, this implies that  $Zhu_{\omega}(V)$ is a semisimple associative algebra.  Note that $Zhu_{\omega}(V)$  has countable dimension. It follows from the proof of Theorem 8.1 of \cite{DLM3} that $Zhu_{\omega}(V)$ is finite dimensional. This completes the proof.
\qed

\section{Rationality of affine vertex operator algebras at admissible levels}
\def\theequation{3.\arabic{equation}}
\setcounter{equation}{0}
In this section, we will show that affine vertex operator algebras at admissible levels are rational with respect to a family of Virasoro elements.
\subsection{Affine vertex operator superalgebras}In this subsection, we recall from \cite{FZ}, \cite{K2}, \cite{LL} some facts about  affine vertex operator superalgebras. Let $\g$ be a finite dimensional simple Lie superalgebra with a nondegenerate even supersymmetric invariant bilinear form $(\cdot|\cdot)$. The affine Lie superalgebra associated to $\g$ is defined on $\hat{\g}=\g\otimes \C[t^{-1}, t]\oplus \C K$ with Lie brackets
\begin{align*}
[x(m), y(n)]&=[x, y](m+n)+(x|y) m\delta_{m+n,0}K,\\
[K, \hat\g]&=0,
\end{align*}
for $x, y\in \g$ and $m,n \in \Z$, where $x(n)$ denotes $x\otimes t^n$.

For a complex number $k$, define the vacuum module of $\hat \g$ by
\begin{align*}
V^{k}(\g)=\Ind_{\g\otimes \C[t]\oplus \C K}^{\hat \g}\C_k,
\end{align*}
where $\C_k=\C\1$ is the $1$-dimensional $\g\otimes \C[t]\oplus \C K$-module such that $\g\otimes \C[t]$ acts as $0$ and $K$ acts as $k$.
\begin{theorem}[\cite{FZ,K2}]\label{affine}
Let $\g$ be a finite dimensional simple Lie superalgebra with a nondegenerate even supersymmetric invariant bilinear form $(\cdot|\cdot)$, $h^\vee$ be the dual Coxeter number of $(\g, (\cdot|\cdot))$ and $k$ be a complex number which is not equal to $-h^\vee$.  Then $V^{k}(\g)$  is a vertex superalgebra. Moreover, if $\{a^i\}$ and $\{b^i\}$ are dual bases of $\g$ with respect to $(\cdot|\cdot)$, i.e., $(b^i|a^j)=\delta_{i,j}$, then
$\omega=\frac{1}{2(k+h^\vee)}\sum_ia^i(-1)b^i(-1)\1$ is a Virasoro element of $V^{k}(\g)$ such that $V^{k}(\g)$ is a $\Z$-graded vertex operator superalgebra.
\end{theorem}

It is well-known that  $V^{k}(\g)$ has a unique maximal proper submodule $J$ \cite{K}. As a result, $V_{k}(\g)=V^{k}(\g)/J$ is a simple vertex superalgebra.
\subsection{Virasoro elements of affine vertex operator algebras} In this subsection, we let $\g$ be a finite dimensional simple Lie algebra. Fix a Cartan subalgebra $\h$ of $\g$ and denote the corresponding root system by $\Delta_{\g}$. We further fix a set of positive roots $\Delta_{\g}^+$ and  denote the  simple roots by $\{\alpha_1,\cdots,\alpha_l\}$.
Then $\g$ has a triangular decomposition $\g=n_-\oplus \h\oplus n_+$, where $n_+$ denotes the subalgebra corresponding to positive roots. Let $\{\alpha_1^\vee,\cdots,\alpha_l^\vee\}$ be the set of simple coroots of $\g$. We then define the fundamental weights $\Lambda_1, \cdots, \Lambda_l$ of $\g$ by  $\langle\Lambda_i, \alpha_j^\vee\rangle=\delta_{i,j}$,  $i,j=1,\cdots, l$.

 Let $\theta$ be a maximal root of $\g$ and $(\cdot|\cdot)$ be the nondegenerate invariant bilinear form of $\g$  such that $(\theta| \theta)=2$. By Theorem \ref{affine},  $V_{k}(\g)$ is a simple vertex algebra. Moreover, $V_{k}(\g)$ has a Virasoro element $\omega$ defined as in Theorem \ref{affine}, which is called the {\em canonical}  Virasoro element of $V_{k}(\g)$. It was proved in \cite{FZ} that $V_{k}(\g)$ equipped with $\omega$ is a $\Z$-graded vertex operator algebra. Furthermore, the Zhu's algebra $Zhu_{\omega}(V^{k}(\g))$  is isomorphic to the universal enveloping algebra $U(\g)$ of $\g$ \cite{FZ}.

For any weak $V_{k}(\g)$-module $M$, we define $\Omega(M)=\{m\in M|(\g\otimes t\C[t])\cdot m=0\}$. Then we have the following result.
\begin{proposition}\label{top}
For any weak $V_{k}(\g)$-module $M$, $\Omega(M)$ is a module of $Zhu_{\omega}(V_{k}(\g))$.
\end{proposition}
\pf Let $M$ be a weak $V_{k}(\g)$-module. We define $$\Omega_{\omega}(M)=\{m\in M|v_{(n)}m=0 \text{ for } v\in V_{k}(\g), n\in \Z \text{ such that } \wt v-n-1<0\}.$$
It was proved in Theorem 5.3 of \cite{DLM3} that $\Omega_{\omega}(M)$ is a module of $Zhu_{\omega}(V_{k}(\g))$.

We next show that $\Omega_{\omega}(M)=\Omega(M)$. First, by the definitions of $\Omega_{\omega}(M)$ and $\Omega(M)$, we have $\Omega_{\omega}(M)\subset \Omega(M)$. We now show that $\Omega_{\omega}(M)\supset \Omega(M)$. Let $w$ be an element in $\Omega(M)$. We need to show that $v_{(n)}w=0$  for  $v\in V_{k}(\g), n\in \Z$  such that  $\wt v-n-1<0$. Note that $v$ has the form $x^1_{(-n_1)}\cdots x^s_{(-n_s)}\1$, where $n_1, \cdots, n_s$ are positive integers. We prove the statement by induction on $s$. When $v=x_{(-m)}\1$, we have $\wt v=m$. If $m-n-1<0$, we have
\begin{align*}
(x_{(-m)}\1)_{(n)}w&=\sum_{j\geq 0}(-1)^j\binom{-m}{j}
(x_{(-m-j)}\1_{(n+j)}-(-1)^m\1_{(n-m-j)}x_{(j)})w\\
&=-(-1)^m\1_{(n-m)}x_{(0)}w=0.
\end{align*}
We now assume that $v_{(n)}w=0$  for any $v=x^1_{(-n_1)}\cdots x^{s-1}_{(-n_{s-1})}\1$ and $n$ such that $\wt v-1<n$, where $n_1, \cdots, n_{s-1}$ are positive integers. We proceed to consider the case that $v=x^1_{(-n_1)}\cdots x^s_{(-n_s)}\1$. Set $u=x^2_{(-n_2)}\cdots x^s_{(-n_s)}\1$, then we have $v=x^1_{(-n_1)}u$. Hence, we have $\wt v=n_1+\wt u$. By the induction assumption, for $n>n_1+\wt u-1$,  we have
\begin{align*}
(x^1_{(-n_1)}u)_{(n)}w&=\sum_{j\geq 0}(-1)^j\binom{-n_1} {j}
(x^1_{(-n_1-j)}u_{(n+j)}-(-1)^{n_1}u_{(n-n_1-j)}x^1_{(j)})w\\
&=-(-1)^{n_1}u_{(n-n_1)}x^1_{(0)}w.
\end{align*}
It is enough to show that $x_{(0)}w\in \Omega(M)$ for any $x\in \g$ and $w\in \Omega(M)$. For any $y\in \g$ and $n\geq 1$, we have $y(n)x(0)w=x(0)y(n)w+[y, x](n)w=0$. This implies that  $x_{(0)}w\in \Omega(M)$ for any $x\in \g$ and $w\in \Omega(M)$. This completes the proof.
\qed

\vskip.25cm
We next show that $V_{k}(\g)$ has a family of Virasoro elements.
Define $$C^{\circ}_{\g, -}=\{h\in \h|\alpha_i(h)\in \Q \text{ and } \alpha_i(h)<0 \text{ for all } 1\leq i\leq l\}.$$
By Proposition \ref{Virasoro}, $V_{k}(\g)$ has a family of Virasoro elements (cf. \cite{AE}, \cite{DLinM}).
\begin{proposition}\label{changeV}
Let $h$ be an element in $C^{\circ}_{\g, -}$, $\sigma$ be a rational number such that $-1<\frac{\sigma}{2}\alpha(h)<0$ for all $\alpha\in \Delta_{\g}^+$. Then $\omega_{\sigma}=\omega-\frac{\sigma}{2}h_{(-2)}\1$ is a Virasoro element of $V_{k}(\g)$. In particular, $V_{k}(\g)$ equipped with $\omega_{\sigma}$ is a vertex operator algebra.
\end{proposition}
\subsection{Affine vertex operator algebras at admissible level}\label{admissiblesec}
 In this subsection, we let $\g$  be a finite dimensional simple Lie algebra, $k$ be an {\em admissible level} of $\hat \g$, i.e., $k=-h^\vee+\frac{p}{q}$, where $p, q$ are positive integers such that $(p, q)=1$ and $$p\geq \left\{
\begin{array}{ll}
h^\vee,&(r^\vee, q)=1,\\
h,& (r^\vee, q)=r^\vee,
\end{array}\right.$$
here $r^\vee$ is the lacing number of $\g$. Our goal in this subsection is to show that $V_k(\g)$ is rational with respect to a family of Virasoro elements.

For $\lambda \in \h^*$, we use $L_{\g}(\lambda)$ to denote the irreducible highest weight module for $\g$ with highest weight $\lambda$ and define
\begin{align*}
V_{\g}(k, \lambda)=\Ind_{\g\otimes \C[t]\oplus \C K}^{\hat \g}L_{\g}(\lambda),
\end{align*}
where $L_{\g}(\lambda)$ is viewed as a module for $\g\otimes \C[t]\oplus \C K$ such that $\g\otimes t\C[t]$ acts as $0$ and $K$ acts as $k$. It is well-known that $V_{\g}(k, \lambda)$ has a unique maximal proper submodule which is denoted by $J(k, \lambda)$ (see \cite{K}). Let $L_{\g}(k, \lambda)$ be the corresponding irreducible quotient module. Then $L_{\g}(k, \lambda)$ is an irreducible weak $V^k(\g)$-module (cf. \cite{FZ}, \cite{LL}).

We now introduce an important category $\mathcal C_{\g,k}$. In the following, for a Lie algebra $\mathfrak{k}$ and a module $M$ of $\mathfrak{k}$, we say that $\mathfrak{k}$ {\em  acts locally nilpotently on} $M$ if for any $m\in M$ there exists a positive integer $s$ such that $x_1\cdots x_s\cdot m=0$ for any $x_1,\cdots ,x_s\in \mathfrak{k}$. We define $\mathcal C_{\g,k}$ to be the subcategory of weak $V_k(\g)$-module category such that $M$ is an object in $\mathcal C_{\g,k}$ if and only if $\g\otimes t\C[t]\oplus n_+$ acts locally nilpotently on  $ M$. We will show that $\mathcal C_{\g,k}$ is semisimple. First, the following result has been established in Theorem 2.20 of \cite{DLM0}.
\begin{theorem}
For $\g=sl_2$, then $\mathcal C_{\g,k}$ is semisimple. Moreover, any object $M$ in $\mathcal C_{\g,k}$ is a direct sum of  modules of the form $L_{\g}(k, \lambda)$, $\lambda\in \h^*$.
\end{theorem}

Recall that $\Omega(M)=\{m\in M|(\g\otimes t\C[t])\cdot m=0\}$ for  a weak $V_k(\g)$-module $M$. Then we have the following result.
\begin{proposition}\label{singular}
Let $\g$ be a finite dimensional simple Lie algebra, $M$ be a weak $V_k(\g)$-module belonging to $\mathcal C_{\g,k}$. Then $\Omega(M)\neq 0$.
\end{proposition}
\pf The proof is similar to that of Theorem 3.7 of \cite{DLM1}. Let $M$ be a weak $V_k(\g)$-module belonging to $\mathcal C_{\g,k}$. For any $m\in M$, set $d(m)=\dim ((\g\otimes \C[t]t)\cdot m)$. Then it enough to show that there exists $0\neq m\in M$ such that $d(m)=0$. Otherwise, assume that for any $0\neq m\in M$, we have $d(m)\neq 0$. Let $m_0$ be a nonzero vector such that $d(m_0)$ is minimal. Our goal is to construct a nonzero vector $w$ such that $d(w)<d(m_0)$.

Since $M$ is a weak $V_k(\g)$-module and $\g$ is finite dimensional,  there exists a positive integer $l$ such that $(\g\otimes t^l\C[t])\cdot m_0=0$. By the assumption that $d(m_0)\neq 0$, there exists a positive integer $j$ such that $(\g\otimes t^j)\cdot m_0\neq 0$ and $(\g\otimes t^n)\cdot m_0=0$ for any $n>j$. In particular, there exists $x\in \g$ such that $x(j)m_0\neq 0$. As a result, by the assumption, there exists a positive integer $r$ such that $x(j)^rm_0\neq 0$ and $x(j)^nm_0= 0$ for any $n>r$.

Set $w=x(j)^rm_0$. We will show that $d(w)<d(m_0)$. First, we show that $$[a,x](n+j)x(j)^sm_0=0$$ holds for any $a\in \g$ and any integers $s\geq 0$, $n\geq 1$. When $s=0$, by the definition of $j$, we have $[a,x](n+j)m_0=0$ holds for any $a\in \g$ and positive integer $n$.  Note that for any positive integer $s$, we have
\begin{align*}
&[a,x](n+j)x(j)^{s+1}m_0\\
&=[a,x](n+j)x(j)x(j)^{s}m_0\\
&=x(j)[a,x](n+j)x(j)^{s}m_0+[[a,x](n+j),x(j)]x(j)^{s}m_0\\
&=x(j)[a,x](n+j)x(j)^{s}m_0+[[a,x],x](n+2j)x(j)^{s}m_0.
\end{align*}
By induction on $s$, $[a,x](n+j)x(j)^sm_0=0$ holds for any $a\in \g$ and any integers $s\geq 0$, $n\geq 1$.

 We next show that if $a(n)m_0=0$ for some $a\in \g$ and positive integer $n$, then $a(n)x(j)^sm_0=0$ holds for any nonnegative integer $s$. Note that for any positive integer $s$, we have
\begin{align*}
&a(n)x(j)^sm_0\\
&=a(n)x(j)x(j)^{s-1}m_0\\
&=x(j)a(n)x(j)^{s-1}m_0+[a(n), x(j)]x(j)^{s-1}m_0\\
&=x(j)a(n)x(j)^{s-1}m_0+[a, x](n+j)x(j)^{s-1}m_0\\
&=x(j)a(n)x(j)^{s-1}m_0.
\end{align*}
By induction on $s$, we have $a(n)x(j)^sm_0=0$ holds for any nonnegative integer $s$. In particular, $a(n)w=0$.

We now consider the linear maps $$\phi_1: \g\otimes\C[t]t\to (\g\otimes\C[t]t)\cdot m_0,~ a(l)\mapsto a(l)m_0,$$$$ \phi_2: \g\otimes\C[t]t\to (\g\otimes\C[t]t)\cdot w,~ a(l)\mapsto a(l)w.$$
By the discussion above, we have $\ker \phi_1\subset \ker \phi_2$. However, $x(j)m_0\neq 0$ and $x(j)w=0$.  This implies that $\ker \phi_1\subsetneqq \ker \phi_2$. As a consequence, $d(w)<d(m_0)$. This is a contradiction. Therefore, we have $\Omega(M)\neq 0$.\qed

We next show that $\Omega(M)$ has a weight vector. The following result has been proved in the proof of Proposition 2.18 and Proposition 2.10 of \cite{DLM0}.
\begin{proposition}\label{a1case}
Let $k$ be an admissible level of $\hat{sl}_2$, $M$ be a weak $V_k(sl_2)$-module such that $\Omega(M)\neq 0$. Assume that $\{e, h, f\}$ be a set of Chevalley generators of $sl_2$ and $w\in \Omega(M)$ is a vector such that $e\cdot w=0$. Then $U(\C h)w$ is finite dimensional and $h$ acts semisimply on $U(\C h)w$.
\end{proposition}

In the general case that $\g$ is a finite dimensional simple Lie algebra,  we have the following result.
\begin{proposition}\label{highest}
Let $\g$ be a finite dimensional simple Lie algebra, $k$ be an admissible level of $\hat \g$, $M$ be a weak $V_k(\g)$-module belonging to $\mathcal C_{\g,k}$. Then  $M$ has a highest weight vector, i.e., there exists $v\in M$ and $\lambda\in \h^*$ such that $(\g\otimes \C[t]t\oplus n_+)\cdot v=0$ and $h\cdot v=\lambda(h)v$ for any $h\in \h$.
\end{proposition}
\pf By Proposition \ref{singular}, we have $\Omega(M)\neq 0$. Thus,  $\Omega(M)$ is a nonzero module of $Zhu_{\omega}(V_k(\g))$. By the definition of $\mathcal C_{\g,k}$, $n_+$ acts locally nilpotently on $M$.  Since $n_+$ is finite dimensional, $U(n_+)m$ is finite dimensional for any $m\in \Omega(M)$, where $U(\cdot)$ denotes the universal enveloping algebra of the Lie algebra. Since $n_+$ is a nilpotent Lie algebra, it follows from Theorem 3.3 of \cite{H} that there exists a vector $w\in U(n_+)m\subset \Omega(M)$ such that $n_+\cdot w=0$. Therefore, $(\g\otimes \C[t]t\oplus n_+)\cdot w=0$.

We next show that  $U(\h)w$ is finite dimensional.  Since $\Omega(M)$ is a nonzero module of $Zhu_{\omega}(V_k(\g))$, by Theorem \ref{Zhu}, $L_{\omega}(\Omega(M))$ is an $\omega$-admissible $V_k(\g)$-module such that $\Omega_{\omega}(L_{\omega}(\Omega(M)))=\Omega(M)$. In particular, we have $(\g\otimes t\C[t])\cdot \Omega(M)=0$.

 For $i=1, \cdots, l$, let $\mathfrak{p}_i$ be the minimal parabolic subalgebra $\C f_i\oplus \h\oplus n_+$, let $\mathfrak{l}_i$ be its Levi subalgebra, let $\mathfrak{m}_i$ be its nilradical. We have
$$\mathfrak{p}_i=\mathfrak{l}_i\oplus \mathfrak{m}_i,\ \ \ \ \mathfrak{l}_i=sl_2^{(i)}\oplus \h_i^{\perp},$$
where $sl_2^{(i)}$ is the copy of $sl_2$ spanned by $e_i, h_i, f_i$, and $\h_i^{\perp}$ is the orthogonal component of $\C h_i$ in $\h$. We have $\mathfrak{m}_i=\oplus_{\alpha\in \Delta_{\g}^+\backslash \{\alpha_i\}}\C x_\alpha$, where $x_{\alpha}$ is a root vector of $\g$ of root $\alpha$.
By Theorem 2.12 of \cite{AFR},  the semi-infinite cohomology space  $H^{\frac{\infty}{2}+0}(\mathfrak{m}_i[t, t^{-1}], L_{\omega}(\Omega(M)))$ is a weak $V_{k_i}(sl_2^{(i)})$-module, where $k_i$ is an admissible level of $\hat{sl_2^{(i)}}$ defined by the formula (28) of \cite{AFR}. The explicit action of $V_{k_i}(sl_2^{(i)})$ on $H^{\frac{\infty}{2}+0}(\mathfrak{m}_i[t, t^{-1}], L_{\omega}(\Omega(M)))$ is given by the formula (14) of \cite{A1}. By the proof of Lemma 4.3 of \cite{A1}, $\Omega(M)^{\mathfrak{m}_i}$ is contained in $H^{\frac{\infty}{2}+0}(\mathfrak{m}_i[t, t^{-1}], L_{\omega}(\Omega(M)))$ (see also Theorem 2.12 of \cite{AFR}). Moreover, by the explicit action of $V_{k_i}(sl_2^{(i)})$ on $H^{\frac{\infty}{2}+0}(\mathfrak{m}_i[t, t^{-1}], L_{\omega}(\Omega(M)))$, we have $(sl_2^{(i)}\otimes\C[t]t)\cdot\Omega(M)^{\mathfrak{m}_i}=0$ (see also Theorem 2.12 of \cite{AFR}). Note that $w\in \Omega(M)^{\mathfrak{m}_i}$ and $e_i\cdot w=0$. It follows from Proposition \ref{a1case} and the formula (14) of \cite{A1} that $U(\C h_i)w$ is finite dimensional.   Therefore, $U(\h)w$ is finite dimensional.

As a result, $U(n_+\oplus \h)w$ is finite dimensional.
By the Lie's Theorem, there exists a vector $v\in U(n_+\oplus \h)w$ such that $n_+\cdot v=0$ and $h\cdot v=\lambda(h)v$ for some $\lambda\in \h^*$ (cf.  Theorem 4.1 of \cite{H}). Since $w\in \Omega(M)$, we have $U(n_+\oplus \h)w\subset \Omega(M)$. Therefore, $v\in \Omega(M)$ is a highest weight vector.\qed

We now classify simple objects in the category $\mathcal C_{\g,k}$.
\begin{proposition}\label{sobject}
Let $\g$ be a finite dimensional simple Lie algebra, $k$ be an admissible level of $\hat \g$. Then any irreducible weak $V_k(\g)$-module belonging to $\mathcal C_{\g,k}$ is of the form $L_{\g}(k, \lambda)$, $\lambda\in \h^*$.
\end{proposition}
\pf Let $M$ be an irreducible weak $V_k(\g)$-module belonging to $\mathcal C_{\g,k}$. By Proposition \ref{highest}, $M$ contains a highest weight vector $v$. Since $M$ is an irreducible weak $V_k(\g)$-module, $M$ is generated by $v$. By Proposition 4.7 of \cite{A1}, $M$ must be of  the form $L_{\g}(k, \lambda)$, $\lambda\in \h^*$.
\qed

We are now ready to prove that  the category $\mathcal C_{\g,k}$ is semisimple.
\begin{theorem}\label{semisimple}
Let $\g$ be a finite dimensional simple Lie algebra, $k$ be an admissible level of $\hat \g$. Then $\mathcal C_{\g,k}$ is semisimple. Moreover, any object $M$ in $\mathcal C_{\g,k}$ is a direct sum of weak $V_k(\g)$-modules of the form $L_{\g}(k, \lambda)$, $\lambda\in \h^*$.
\end{theorem}
\pf Let $M$ be a weak $V_k(\g)$-module belonging to $\mathcal C_{\g,k}$. By Proposition \ref{highest}, $M$ contains a highest weight vector $v$. By Proposition 4.7 of \cite{A1}, the weak $V_k(\g)$-submodule $\langle v \rangle$ of $M$ generated by $v$ is of the form  $L_{\g}(k, \lambda)$, $\lambda\in \h^*$. We now let $W$ be the sum of irreducible weak $V_k(\g)$-submodules of $M$. Then $W$ is a direct sum of weak $V_k(\g)$-modules of the form $L_{\g}(k, \lambda)$, $\lambda\in \h^*$.

 Our goal is to show that $W=M$. Otherwise, $W$ is a proper submodule of $M$. As a result, $M/W$ is a weak $V_k(\g)$-module belonging to $\mathcal C_{\g,k}$. By Proposition \ref{highest}, $M/W$ has a highest weight vector $\bar{w}$. Let $w$ be a preimage of $\bar w$, then we have $n_+\cdot w\subset W$. Since $n_+$ is finitely generated, there exist submodules $L_{\g}(k, \lambda_1), \cdots, L_{\g}(k, \lambda_s)$ of $W$ such that $n_+w\subset  L_{\g}(k, \lambda_1)\oplus \cdots\oplus L_{\g}(k, \lambda_s)$. Therefore, it follows from  Proposition 4.7 of \cite{A1} that the submodule of $M$ generated by $w$ and $ L_{\g}(k, \lambda_1)\oplus \cdots\oplus L_{\g}(k, \lambda_s)$ is completely reducible. In particular, the submodule of  $M$ generated by $w$ is a direct sum of irreducible weak $V_k(\g)$-submodules of the form $L_{\g}(k, \lambda)$, $\lambda\in \h^*$. This is a contradiction. Thus, $W=M$, this completes the proof.
\qed

We are ready to prove the main result in this section.
\begin{theorem}\label{rationalaff}
Let $h$ be an element in $C^{\circ}_{\g, -}$, $\sigma$ be a rational number such that $-1<\frac{\sigma}{2}\alpha(h)<0$ for all $\alpha\in \Delta_{\g}^+$, and $k$ be an admissible level of $\hat \g$. Then $V_k(\g)$ is $\omega_{\sigma}$-rational. In particular, the Zhu's algebra $Zhu_{\omega_{\sigma}}(V_k(\g))$ is semisimple.
\end{theorem}
\pf Let $M$ be an $\omega_{\sigma}$-admissible $V_k(\g)$-module. Note that for any root $\alpha\in \Delta_{\g}$ and the corresponding root vector  $x_\alpha$, we have
\begin{align*}
(\omega_{\sigma})_{(1)}x_\alpha
&=\omega_{(1)}x_\alpha+\frac{\sigma}{2}h_{(0)}x_\alpha\\
&=(1+\frac{\sigma}{2}\alpha(h))x_\alpha.
\end{align*}
By the assumption,  $-1<\frac{\sigma}{2}\alpha(h)<0$ for all $\alpha\in \Delta_{\g}^+$, this implies that $M$ belongs to $\mathcal{C}_{\g, k}$. By Theorem \ref{semisimple}, $M$ is a direct sum of  weak $V_k(\g)$-modules of the form $L_{\g}(k, \lambda)$, $\lambda\in \h^*$.

We now prove that weak $V_k(\g)$-modules of the form $L_{\g}(k, \lambda)$ are $\omega_{\sigma}$-admissible. Let $w$ be a highest weight vector of  $L_{\g}(\lambda)$. By the definition of $L_{\g}(k, \lambda)$ and the PBW Theorem, $L_{\g}(k, \lambda)$ is spanned by the vectors of the form $x^1(-n_1)\cdots x^s(-n_s)y^1(0)\cdots y^t(0) w$, where $x^1, \cdots, x^s\in \g$, $n_1, \cdots, n_s\in \Z_{> 0}$ and $y^1, \cdots, y^t\in n_{-}$. It follows that $L_{\g}(k, \lambda)$ is an $\omega_{\sigma}$-ordinary module of $V_k(\g)$. By Proposition \ref{ordi}, $L_{\g}(k, \lambda)$ is an $\omega_{\sigma}$-admissible module of $V_k(\g)$. In particular, $M$ is a direct sum of  irreducible $\omega_{\sigma}$-admissible  $V_k(\g)$-modules. This completes the proof.
\qed
\section{Extensions of vertex operator algebras}
We will show that certain affine vertex operator superalgebras and minimal $W$-algebras are rational with respect to a family of Virasoro elements in next section. To obtain these results, we prove that certain categories  of weak modules of  extensions of vertex operator algebras are semisimple.

 In this section, we assume that  $(U,Y(\cdot, z),\1,\omega)$ and $(V,Y(\cdot, z),\1,\omega)$ are vertex operator superalgebras satisfying the following conditions:\\
(1) $(V,Y(\cdot, z),\1,\omega)$ is a simple vertex operator algebra;\\
(2) $(V,Y(\cdot, z),\1,\omega)$ is a subalgebra of $(U,Y(\cdot, z),\1,\omega)$ having the same Virasoro element. \\
(3) There exists an irreducible $\omega$-ordinary $V$-module $V^1$ such that $U$ viewed as an $\omega$-ordinary $V$-module has the decomposition $U=V\oplus V^1$. Moreover, $V^1\cdot V^1=\{u_{(n)}w|u, w\in V^1, n\in \Z\}$ viewed as an $\omega$-ordinary $V$-module is isomorphic to $V$.

Let $\mathcal C$ be a subcategory of weak $V$-module category such that\\
(1) $\mathcal C$ is semisimple.\\
(2) For any weak $V$-module $M$ in $\mathcal C$, all the $V$-submodules of $M$ are also in $\mathcal C$.\\
 We then define $\mathcal D$ to be the subcategory of weak $U$-module category such that $M\in \mathcal D$ if and only if $M$ viewed as a weak $V$-module belongs to $\mathcal C$. In this section, our goal is to prove that $\mathcal D$ is semisimple.

Let $M$ be an object in $\mathcal D$. Then $M$ viewed as a weak $V$-module belongs to $\mathcal C$. Thus, $M$ viewed as a weak $V$-module  is completely reducible, i.e., $M=\oplus_{i\in I}M^i$, where $M^i~(i\in I)$ are irreducible weak $V$-modules. In the following, for any $A\subset U$ and $B\subset M$, we use $A\cdot B$ to denote the set $\{u_{(n)}w|u\in A, w\in B, n\in \Z\}$. Then we have
\begin{lemma}\label{scur}
For any $\lambda\in I$, $V^1\cdot M^{\lambda}$ is an irreducible weak $V$-module.
\end{lemma}
\pf Let $u, v\in U$, $p, q\in \Z$ and $w\in M^{\lambda}$. Since $M$ is a weak $U$-module, it follows from the proof of Proposition 4.5.7 of \cite{LL} that $u_{(p)}v_{(q)}w$ can be expressed as a linear combination of elements of the form $t_{(r)}w$, and $t$ can be chosen to be of the form $u_{(s)}v$ for $s\in \Z$. This implies that $V^1\cdot M^{\lambda}$ is a weak $V$-module.

We next show that $V^1\cdot M^{\lambda}$ is  irreducible. Otherwise, since $V^1\cdot M^{\lambda}$ is an object in $\mathcal{C}$, $V^1\cdot M^{\lambda}$ is a direct sum of irreducible weak $V$-modules, i.e., $V^1\cdot M^{\lambda}= W^1\oplus\cdots\oplus W^k$, where $W^i, 1\leq i\leq k,$ are irreducible weak $V$-modules and $k\geq 2$.  Then we have $$ V^1\cdot(V^1\cdot M^{\lambda})=V^1\cdot W^1+\cdots+ V^1\cdot W^k.$$

We now show that $V^1\cdot W^i\neq 0$, $1\leq i\leq k$. Otherwise, assume that $V^1\cdot W^i= 0$ for some $i$. Then we have $V^1\cdot (V^1\cdot W^i)=0$. However, $(V^1\cdot V^1)\cdot W^i=V^1\cdot (V^1\cdot W^i)=0$ and $V^1\cdot V^1=V$, this is a contradiction. Hence, $V^1\cdot W^i\neq 0$, $1\leq i\leq k$.

We now show that $$V^1\cdot(V^1\cdot M^{\lambda})=V^1\cdot W^1\oplus\cdots\oplus V^1\cdot W^k.$$ Otherwise, there exists $i$ such that $V^1\cdot W^i\cap \sum_{j\neq i}V^1\cdot W^j\neq 0$. By the similar argument as above, we have $V^1\cdot (V^1\cdot W^i\cap \sum_{j\neq i}V^1\cdot W^j)\neq 0$.  Hence, there exist $u\in V^1$, $w\in V^1\cdot W^i\cap \sum_{j\neq i}V^1\cdot W^j$ and $n\in \Z$ such that $u_{(n)}w\neq 0$. Since $w\in V^1\cdot W^i\cap \sum_{j\neq i}V^1\cdot W^j$, there exist $s_i\in \Z~(1\leq i\leq k)$, $b^i\in V^1~(1\leq i\leq k)$ and $w^i\in W^i~(1\leq i\leq k)$ such that $$w=b^i_{(s_i)}w^i=\sum_{j\neq i} b^j_{(s_j)}w^j.$$ Therefore, we have$$0\neq u_{(n)}w=u_{(n)}b^i_{(s_i)}w^i=\sum_{j\neq i} u_{(n)}b^j_{(s_j)}w^j.$$ Note that $u_{(n)}b^l_{(s_l)}w^l\in W^l$ for $1\leq l\leq k$. However, we have $V^1\cdot M^{\lambda}= W^1\oplus\cdots\oplus W^k$. This is a contradiction. Then we have $ V^1\cdot(V^1\cdot M^{\lambda})=V^1\cdot W^1\oplus\cdots\oplus V^1\cdot W^k$ and $k\geq 2$.
However, $V^1\cdot(V^1\cdot M^{\lambda})\subset (V^1\cdot V^1)\cdot M^{\lambda}=M^{\lambda}$. This is a contradiction. Then $V^1\cdot M^{\lambda}$ is  irreducible.
\qed

\vskip.25cm

We are now ready to prove the main result in this section.
\begin{theorem}\label{scextension}
Let $\mathcal C$ and $\mathcal D$ be categories defined as above. Then the category $\mathcal D$ is semisimple.
\end{theorem}
\pf Let $M$ be an object in $\mathcal D$. Then $M$ viewed as a weak $V$-module belongs to $\mathcal C$. Thus, $M$ viewed as a weak $V$-module  is completely reducible, i.e., $M=\oplus_{i\in I}M^i$, where $M^i~(i\in I)$ are irreducible weak $V$-modules. For any $M^i$, we let $\langle M^i\rangle$ be the weak $U$-submodule of $M$ generated by $M^i$. Then $\langle M^i\rangle=M^i+V^1\cdot M^i$. By the discussion above, it is enough to show that $\langle M^i\rangle$ viewed as a weak $U$-module is completely reducible. By Lemma \ref{scur}, we have $\langle M^i\rangle=M^i$ or $\langle M^i\rangle=M^i\oplus V^1\cdot M^i$. If $\langle M^i\rangle=M^i$, then $\langle M^i\rangle$ is an irreducible weak $U$-module. We then consider the case that $\langle M^i\rangle=M^i\oplus V^1\cdot M^i$. Therefore, for any $w\in \langle M^i\rangle$, there exist unique elements $w^0\in M^i$ and $w^1\in V^1\cdot M^i$ such that $w=w^0+w^1$. For a nontrivial proper weak $U$-submodule $W$ of $\langle M^i\rangle$, we define two linear maps $P_0, P_1$ as follows:
\begin{align*}
P_0: W\to M^i,~
w\mapsto w^0,
~~~~~~~~~~~~~~~~\ \ \ \ \ \
P_1:W\to  V^1\cdot M^i,
~w\mapsto w^1.
\end{align*}
Since $M^i$ and  $V^1\cdot M^i$ are weak $V$-modules, $P_0, P_1$ are $V$-module homomorphisms. In particular, $\ker P_0$ and $\ker P_1$ are weak $V$-submodules of $V^1\cdot M^i$ and $M^i$, respectively. We now show that
$\ker P_0=0$. Otherwise, $\ker P_0\neq 0$. Note that $\ker P_0=W\cap V^1\cdot M^i$ Since $V^1\cdot M^i$ is an irreducible weak $V$-module, we have  $\ker P_0=W\cap V^1\cdot M^i=V^1\cdot M^i$. This implies that $V^1\cdot M^i\subset W$. Since $W$ is a weak $U$-module, we have $M^i\subset W$. Hence, $W=\langle M^i\rangle$. This is a contradiction, and we have $\ker P_0=0$.  Similarly, $\ker P_1=0$.

We next show that $\im P_0\neq 0$. Otherwise, $\im P_0= 0$, this implies that $W\subset V^1\cdot M^i$. By Lemma \ref{scur}, $V^1\cdot M^i$ is an irreducible weak $V$-module, this implies that $W=V^1\cdot M^i$. As a result, we have $M^i\subset W$. However, $W$ is a proper weak $U$-submodule of $\langle M^i\rangle$.  This is a contradiction. Therefore, $\im P_0\neq 0$. Similarly, $\im P_1\neq 0$. Since $M^i$ and  $V^1\cdot M^i$ are irreducible weak $V$-modules, we have $\im P_0=M^i$ and $\im P_1= V^1\cdot M^i$. As a consequence, $P_0, P_1$ induce two $V$-module isomorphisms $P_0: W\to M^i$ and $P_1: W\to V^1\cdot M^i$.

Note that we have $M^i\cong V^1\cdot M^i$ or $M^i\ncong V^1\cdot M^i$. If $M^i\ncong V^1\cdot M^i$, by the discussion above, $\langle M^i\rangle$ has no nontrivial proper weak $U$-submodules. Thus, $\langle M^i\rangle$ is an irreducible weak $U$-module.

We then consider the case that  $M^i\cong V^1\cdot M^i$. By the discussion above, if $W$ is  a nontrivial proper weak $U$-submodule of $\langle M^i\rangle$, then $P_1P_0^{-1}: M^i\to V^1\cdot M^i$ is a $V$-module isomorphism. Note that for any $w\in W$, we have $w=w^0+w^1$, $w^0\in M^i$, $w^1\in V^1\cdot M^i$. In particular, $w=P_0(w)+P_1(w)=P_0(w)+P_1(P_0^{-1}P_0(w))$. Since $P_0:W\to M^i$ is a $V$-module isomorphism, we have $$W=\{w^0+P_1P_0^{-1}(w^0)|w^0\in M^i\}.$$
Since $W$ is a weak $U$-module, for any $u\in V^1$ and $w^0\in M^i$, we have
$Y(u, z)(w^0+P_1P_0^{-1}(w^0))\in W$. This implies that $Y(u, z)w^0+Y(u, z)P_1P_0^{-1}(w^0)\in W$. Therefore, we have
\begin{align*}
Y(u, z)w^0=P_1P_0^{-1}Y(u, z)P_1P_0^{-1}(w^0).
 \end{align*}
Set $\tilde{W}=\{w^0-P_1P_0^{-1}(w^0)|w^0\in M^i\}.$ Since $P_1P_0^{-1}$  is a $V$-module isomorphism, $\tilde{W}$ is a weak $V$-submodule of $\langle M^i\rangle$. Furthermore, for any $u\in V^1$ and $w^0\in M^i$, we have $Y(u, z)w^0=P_1P_0^{-1}Y(u, z)P_1P_0^{-1}(w^0)$. This implies $\tilde{W}$ is a weak $U$-submodule of $\langle M^i\rangle$.

By the definition of $\tilde{W}$, we have $\langle M^i\rangle=W+\tilde{W}$. We next show that $W\cap \tilde{W}=0$. Note that if $w^0$ is an element of $ M^i$ such that $w^0+P_1P_0^{-1}(w^0)\in W\cap \tilde{W}$, we have $P_1P_0^{-1}(w^0)=-P_1P_0^{-1}(w^0)$. This forces that $P_1P_0^{-1}(w^0)=0$ and $w^0=0$. Therefore, $W\cap \tilde{W}=0$ and we have $\langle M^i\rangle=W\oplus\tilde{W}$. Note that $W$ and $ \tilde{W}$ are nontrivial proper weak $U$-submodules of $\langle M^i\rangle$. By the discussion above, $W$ and $ \tilde{W}$ viewed as  weak $V$-modules are isomorphic to $M^i$. Thus, $W$ and $ \tilde{W}$ are irreducible weak $U$-module. This completes the proof.
\qed

\begin{remark}
When $V$ is a regular $\Z$-graded vertex operator algebra, $V^1$ is a simple current module of $V$, and the results in this section have been obtained in \cite{L3}. In general, the category $\mathcal C$ may not be a tensor category, and we need to modify the arguments in \cite{L3}.
\end{remark}

\section{Rationality of certain affine vertex operator superalgebras}\label{affineVOSA}
In this section, we will show that certain affine vertex operator superalgebras  are rational with respect to a family of  Virasoro elements.
\subsection{Conformal embeddings}
 Let $\g=\g_{\bar 0}\oplus \g_{\bar 1}$ be a basic classical simple Lie superalgebra, $(\cdot|\cdot)$ be the non-degenerate invariant supersymmetric bilinear form defined in Table 6.1 of \cite{KW4},  $h^\vee$ be the dual Coxeter number of $(\g, (\cdot|\cdot))$ and $k$ be a complex number which is not equal to $-h^\vee$.  Then $V_{k}(\g)$  is a simple vertex superalgebra. Moreover, $V_{k}(\g)$ has a canonical Virasoro element $\omega_{\g}$ defined as in Theorem \ref{affine}. We use $\mathcal{V}_{k}(\g_{\bar 0})$ to denote the vertex subalgebra of $V_{k}(\g)$  generated by $x(-1)\1$, $x\in \g_{\bar 0}$. Then there exists a complex number $u(k)$ such that $\mathcal{V}_{k}(\g_{\bar 0})$ is a quotient of $V^{u(k)}(\g_{\bar 0})$. Let $\omega_{\g_{\bar 0}}$ be the canonical Virasoro element of $V^{u(k)}(\g_{\bar 0})$ defined as in Theorem \ref{affine}.  We say that $\mathcal{V}_{k}(\g_{\bar 0})$ is conformally embedded in $V_{k}(\g)$ if $\omega_{\g_{\bar 0}}=\omega_{\g}$. The following results about conformal embeddings were proved in Proposition  4.1 of \cite{AMPP}.
 \begin{proposition}\label{embedding}
 (1) The vertex operator algebra $V_1(sl_2)\otimes V_{-4/3}(G_2)$ is conformally embedded in the vertex operator superalgebra $V_{-4/3}(G(3))$. Moreover, $V_{-4/3}(G(3))$ viewed as an $\omega_{G(3)}$-ordinary $V_1(sl_2)\otimes V_{-4/3}(G_2)$-module has  the following decomposition $$ V_{-4/3}(G(3))=V_1(sl_2)\otimes V_{-4/3}(G_2)\oplus L_{sl_2}(1, \Lambda_1)\otimes L_{G_2}(-4/3, \Lambda_1).$$
 (2) The vertex operator algebra $V_1(sl_2)\otimes V_{-3/2}(so(7))$ is conformally embedded in the vertex operator superalgebra $V_{-3/2}(F(4))$. Moreover, $V_{-3/2}(F(4))$ viewed as an $\omega_{F(4)}$-ordinary $V_1(sl_2)\otimes V_{-3/2}(so(7))$-module has  the following decomposition $$ V_{-3/2}(F(4))=V_1(sl_2)\otimes V_{-3/2}(so(7))\oplus L_{sl_2}(1, \Lambda_1)\otimes L_{so(7)}(-3/2, \Lambda_3).$$
 (3) For $m\neq n$, the vertex operator algebra $V_1(so(2m+1))\otimes V_{-1/2}(sp(2n))$ is conformally embedded in the vertex operator superalgebra $V_{1}(B(m,n))$. Moreover, $V_{1}(B(m,n))$ viewed as an $\omega_{B(m,n)}$-ordinary $V_1(so(2m+1))\otimes V_{-1/2}(sp(2n))$-module has  the following decomposition $$V_{1}(B(m,n))=V_1(so(2m+1))\otimes V_{-1/2}(sp(2n))\oplus L_{so(2m+1)}(1, \Lambda_1)\otimes L_{sp(2n)}(-1/2, \Lambda_1).$$
 (4) The vertex operator algebra $V_{-(2n+3)/4}(sp(2n))$ is conformally embedded in the vertex operator superalgebra $V_{-(2n+3)/2}(B(0,n))$. Moreover, $V_{-(2n+3)/2}(B(0,n))$ viewed as an $\omega_{B(0,n)}$-ordinary $V_{-(2n+3)/4}(sp(2n))$-module has  the following decomposition $$V_{-(2n+3)/2}(B(0,n))=V_{-(2n+3)/4}(sp(2n))\oplus  L_{sp(2n)}(-(2n+3)/4, \Lambda_1).$$
 (5) For $m\neq n+1$, the vertex operator algebra $V_1(so(2m))\otimes V_{-1/2}(sp(2n))$ is conformally embedded in the vertex operator superalgebra $V_{1}(D(m,n))$. Moreover, $V_{1}(D(m,n))$ viewed as an $\omega_{D(m,n)}$-ordinary $V_1(so(2m))\otimes V_{-1/2}(sp(2n))$-module has  the following decomposition $$V_{1}(D(m,n))=V_1(so(2m))\otimes V_{-1/2}(sp(2n))\oplus L_{so(2m)}(1, \Lambda_1)\otimes L_{sp(2n)}(-1/2, \Lambda_1).$$
 \end{proposition}

\subsection{Rationality of affine vertex operator superalgebras at conformal level} Our goal in this subsection is to show that the vertex superalgebras in Proposition \ref{embedding} are rational with respect to a family of Virasoro elements. We will prove the rationality by using Theorem \ref{scextension}. First, let $\mathcal{C}_{sp(2n),-(2n+3)/4}$ be the subcategory of weak $V_{-(2n+3)/4}(sp(2n))$-module category defined in subsection \ref{admissiblesec}. We define $\mathcal{D}_{B(0, n),-(2n+3)/2}$ to be the subcategory of weak $V_{-(2n+3)/2}(B(0, n))$-module category such that $M$ is an object in $\mathcal{D}_{B(0, n),-(2n+3)/2}$ if and only if $M$ viewed as a weak $V_{-(2n+3)/4}(sp(2n))$-module belongs to  $\mathcal{C}_{sp(2n),-(2n+3)/4}$. Then we have the following result.
 \begin{proposition}\label{D-B(0,n)}
 The category $\mathcal{D}_{B(0, n),-(2n+3)/2}$ is semisimple.
 \end{proposition}
\pf By Proposition \ref{embedding}, the vertex operator algebra $V_{-(2n+3)/4}(sp(2n))$ is conformally embedded in the vertex operator superalgebra $V_{-(2n+3)/2}(B(0,n))$. Moreover, $V_{-(2n+3)/2}(B(0,n))$ viewed as an $\omega_{B(0,n)}$-ordinary $V_{-(2n+3)/4}(sp(2n))$-module has  the following decomposition $$V_{-(2n+3)/2}(B(0,n))=V_{-(2n+3)/4}(sp(2n))\oplus  L_{sp(2n)}(-(2n+3)/4, \Lambda_1).$$
Furthermore, by Theorem 2.3 of \cite{AMPP}, $V_{-(2n+3)/4}(sp(2n))$ and $L_{sp(2n)}(-(2n+3)/4, \Lambda_1)$ are the even and odd parts of $V_{-(2n+3)/2}(B(0,n))$, respectively. This implies that $$L_{sp(2n)}(-(2n+3)/4, \Lambda_1)\cdot L_{sp(2n)}(-(2n+3)/4, \Lambda_1)=V_{-(2n+3)/4}(sp(2n)).$$

Since $-(2n+3)/4$ is an admissible level of $\hat{sp(2n)}$, it follows from Theorem \ref{semisimple} that $\mathcal{C}_{sp(2n),-(2n+3)/4}$ is semisimple. Moreover, if $M$ is a weak $V_{-(2n+3)/4}(sp(2n))$-module belonging to  $\mathcal{C}_{sp(2n),-(2n+3)/4}$, then any weak $V_{-(2n+3)/4}(sp(2n))$-submodule of $M$ belongs to  $\mathcal{C}_{sp(2n),-(2n+3)/4}$. It follows from Theorem \ref{scextension} that the category $\mathcal{D}_{B(0, n),-(2n+3)/2}$ is semisimple.
\qed

We now show that the vertex superalgebra $V_{-(2n+3)/2}(B(0, n))$ is rational with respect to a family of Virasoro elements. First, we show that $V_{-(2n+3)/2}(B(0, n))$  has a family of Virasoro elements.
\begin{proposition}\label{v-B(0,n)}
Let $h$ be an element in $C^{\circ}_{sp(2n), -}$, $\sigma$ be a rational number such that $-1<\frac{\sigma}{2}\alpha(h)<0$ for all $\alpha\in \Delta_{sp(2n)}^+$, $\omega_{sp(2n)}$ be the canonical Virasoro element of $V_{-(2n+3)/4}(sp(2n))$ defined as in Theorem \ref{affine}. Set $\omega_{sp(2n), \sigma}=\omega_{sp(2n)}-\frac{\sigma}{2}h_{(-2)}\1$. Then $\omega_{sp(2n), \sigma}$ is a Virasoro element of $V_{-(2n+3)/2}(B(0, n))$.
\end{proposition}
\pf By Proposition \ref{changeV}, $\omega_{sp(2n), \sigma}$ is a Virasoro element of $V_{-(2n+3)/4}(sp(2n))$. By Proposition \ref{embedding}, $V_{-(2n+3)/2}(B(0,n))=V_{-(2n+3)/4}(sp(2n))\oplus  L_{sp(2n)}(-(2n+3)/4, \Lambda_1).$ Moreover, it has been proved in the proof of Theorem \ref{rationalaff} that $L_{sp(2n)}(-(2n+3)/4, \Lambda_1)$ is an $\omega_{sp(2n), \sigma}$-ordinary module of $V_{-(2n+3)/4}(sp(2n))$. This implies that $\omega_{sp(2n), \sigma}$ is a Virasoro element of $V_{-(2n+3)/2}(B(0, n))$.
\qed

We are now ready to prove the rationality of $V_{-(2n+3)/2}(B(0, n))$.
\begin{theorem}\label{r-B(0,n)}
 The vertex superalgebra $V_{-(2n+3)/2}(B(0, n))$ is $\omega_{sp(2n), \sigma}$-rational.
\end{theorem}
\pf Let $M$ be an $\omega_{sp(2n), \sigma}$-admissible $V_{-(2n+3)/2}(B(0, n))$-module. By the argument in the proof of Theorem \ref{rationalaff}, $M$ viewed as a weak $V_{-(2n+3)/4}(sp(2n))$-module belongs to $\mathcal{C}_{sp(2n),-(2n+3)/4}$. Therefore, $M$ belongs to  $\mathcal{D}_{B(0, n),-(2n+3)/2}$. By Proposition \ref{D-B(0,n)}, $M$ is a direct sum of irreducible weak $V_{-(2n+3)/2}(B(0, n))$-modules.

We now show that $M$ is a direct sum of irreducible $\omega_{sp(2n), \sigma}$-admissible  modules of $V_{-(2n+3)/2}(B(0, n))$. By Proposition \ref{sobject}, irreducible weak $V_{-(2n+3)/4}(sp(2n))$-modules belonging to $\mathcal{C}_{sp(2n),-(2n+3)/4}$ are of the form $L_{sp(2n)}(-(2n+3)/4, \lambda)$. It follows from the Main Theorem of \cite{A1} that the number of irreducible weak $V_{-(2n+3)/4}(sp(2n))$-modules belonging to $\mathcal{C}_{sp(2n),-(2n+3)/4}$ is finite. It has been proved in the proof of Theorem \ref{rationalaff} that $L_{sp(2n)}(-(2n+3)/4, \lambda)$ is an $\omega_{sp(2n), \sigma}$-ordinary module of $V_{-(2n+3)/4}(sp(2n))$. This implies that $(\omega_{sp(2n), \sigma})_{(1)}$ acts semisimply on $M$. Moreover, $M_{\mu}=\{m\in M|(\omega_{sp(2n), \sigma})_{(1)}m=\mu m\}$ is equal to $0$ if the real part of $\mu$ is sufficiently small.   By the similar argument as in the proof of Lemma 3.4 of \cite{DLM3}, each irreducible weak $V_{-(2n+3)/2}(B(0, n))$-modules is $\omega_{sp(2n), \sigma}$-admissible. This completes the proof.
\qed

\vskip.25cm
Let $\g$ be a basic classical simple Lie superalgebra listed in TABLE \ref{tab1}. Then we define the Lie subalgebras $\g_1$ and $\g_2$ of $\g$ as in TABLE \ref{tab1}. By Proposition \ref{embedding}, the vertex operator algebra $V_{k_1}(\g_1)\otimes V_{k_2}(\g_2)$ is conformally embedded in $V_k(\g)$, where $k, k_1, k_2$ are the complex numbers listed in TABLE \ref{tab1}. Let $\mathcal{C}_{\g_1, \g_2}$ be the subcategory of weak $V_{k_1}(\g_1)\otimes V_{k_2}(\g_2)$-module category such that $M$ is an object in $\mathcal{C}_{\g_1, \g_2}$ if and only if $M$ viewed as a weak $V_{k_2}(\g_2)$-module belongs to $\mathcal{C}_{\g_2, k_2}$.

\begin{table}[h]
\centering
\begin{tabular}{|c c c c c c|}
\hline
$\g$ & $\g_1$& $\g_2$& $k_1$& $k_2$& $k$\\
\hline
$G(3)$ &  $sl_2$& $G_2$& $1$& $-4/3$& $-4/3$\\
\hline
$F(4)$ &  $sl_2$& $so(7)$& $1$& $-3/2$& $-3/2$\\
\hline
$B(m,n), m\neq n$ &  $so(2m+1)$& $sp(2n)$& $1$& $-1/2$& $1$\\
\hline
$D(m,n), m\neq n+1$ &  $so(2m)$& $sp(2n)$& $1$& $-1/2$& $1$\\
 \hline
\end{tabular}
\caption{}
\label{tab1}
\end{table}
\begin{proposition}\label{semisimpleten1}
Let $\g_1$, $\g_2$ and $\mathcal{C}_{\g_1, \g_2}$ be as above. Then the category $\mathcal{C}_{\g_1, \g_2}$ is semisimple.
\end{proposition}
\pf The proof is similar to that of Proposition 2.7 of \cite{DMZ}. Let $M$ be a weak $V_{k_1}(\g_1)\otimes V_{k_2}(\g_2)$-module. For any $w\in M$, let $\langle w\rangle$ be the weak $V_{k_1}(\g_1)\otimes V_{k_2}(\g_2)$-submodule of $M$ generated by $w$. It is enough to prove that $\langle w\rangle$ is completely reducible. By the formula (4.7.4) of \cite{FHL}, $\langle w\rangle$  is spanned by elements of the form
$$(v_{1,1}\otimes \1)_{(m_1)}\cdots (v_{1,s}\otimes \1)_{(m_s)}(\1 \otimes v_{2,1})_{(n_1)}\cdots (\1 \otimes v_{2,l})_{(n_l)}w,$$
where $v_{1,i}\in V_{k_1}(\g_1)$ and $v_{2,j}\in V_{k_2}(\g_2)$. Set $$M_1={\rm span}\{(v_{1,1}\otimes \1)_{(m_1)}\cdots (v_{1,s}\otimes \1)_{(m_s)}w|v_{1,i}\in V_{k_1}(\g_1), m_i\in \Z\}$$
and $$M_2={\rm span}\{(\1 \otimes v_{2,1})_{(n_1)}\cdots (\1 \otimes v_{2,l})_{(n_l)}w|v_{2,j}\in V_{k_2}(\g_2), n_j\in \Z\}.$$
Then $M_1$ and $M_2$ are weak modules of $V_{k_1}(\g_1)$ and $V_{k_2}(\g_2)$, respectively. Since $V_{k_1}(\g_1)$ is a regular vertex operator algebra, then  $M_1$ is completely reducible. By the assumption, $M_2$ is a weak $V_{k_2}(\g_2)$-module belonging to $\mathcal{C}_{\g_2, k_2}$. Since $k_2$ is an admissible level of $\hat{\g_2}$,  $\mathcal{C}_{\g_2, k_2}$ is semisimple by Theorem \ref{semisimple}. In particular, $M_2$ is completely reducible. By the  formula (4.7.21) of \cite{FHL}, there exists a weak  $V_{k_1}(\g_1)\otimes V_{k_2}(\g_2)$-module epimorphism $\phi: M_1\otimes M_2\to \langle w\rangle$. This implies that $\langle w\rangle$ is completely reducible.
\qed

Let $\g$ be a  basic classical simple Lie superalgebra listed in TABLE \ref{tab1}, $k, k_1, k_2$ be the complex numbers listed in TABLE \ref{tab1}. Then we define $\mathcal{D}_{\g, k}$ to be the subcategory of weak $V_{k}(\g)$-module category such that $M$ is an object in $\mathcal{D}_{\g, k}$ if and only if $M$ viewed as a weak $V_{k_1}(\g_1)\otimes V_{k_2}(\g_2)$-module belongs to $\mathcal{C}_{\g_1, \g_2}$.
\begin{proposition}\label{D-AS}
Let $\g$ and $k$ be as above. Then the category $\mathcal{D}_{\g, k}$ is semisimple.
\end{proposition}
\pf  The proof is similar to that of Proposition \ref{D-B(0,n)}. By Proposition \ref{embedding}, the vertex operator algebra $V_{k_1}(\g_1)\otimes V_{k_2}(\g_2)$ is conformally embedded in the vertex operator superalgebra $V_{k}(\g)$. Moreover, by Theorem 2.3 of \cite{AMPP}, $V_{k}(\g)$ viewed as an $\omega_{\g}$-ordinary $V_{k_1}(\g_1)\otimes V_{k_2}(\g_2)$-module has  the decomposition $$V_{k}(\g)=V_{k_1}(\g_1)\otimes V_{k_2}(\g_2)\oplus  V^1,$$ where $V^1$ is the odd part of $V_{k}(\g)$. In particular, $V^1$ is an irreducible $\omega_{\g}$-ordinary $V_{k_1}(\g_1)\otimes V_{k_2}(\g_2)$-module such that $V^1\cdot V^1=V_{k_1}(\g_1)\otimes V_{k_2}(\g_2)$.
Moreover, if $M$ is a weak $V_{k_1}(\g_1)\otimes V_{k_2}(\g_2)$-module belonging to  $\mathcal{C}_{\g_1,\g_2}$, then any weak $V_{k_1}(\g_1)\otimes V_{k_2}(\g_2)$-submodule of $M$ belongs to  $\mathcal{C}_{\g_1,\g_2}$. It follows from Proposition \ref{semisimpleten1} and Theorem \ref{scextension} that the category $\mathcal{D}_{\g, k}$ is semisimple.
\qed

We now show that the vertex superalgebra $V_{k}(\g)$ is rational with respect to a family of Virasoro elements. First, we show that $V_{k}(\g)$  has a family of Virasoro elements.
\begin{proposition}\label{vir-aff-ten}
Let $h$ be an element in $C^{\circ}_{\g_2, -}$, $\sigma$ be a rational number such that $-1<\frac{\sigma}{2}\alpha(h)<0$ for all $\alpha\in \Delta_{\g_2}^+$, $\omega_{\g_1}$ and $\omega_{\g_2}$ be the canonical Virasoro elements of $V_{k_1}(\g_1)$ and $V_{k_2}(\g_2)$, respectively. Set $\omega_{\g, \sigma}=\omega_{\g_1}\otimes \1+1\otimes \omega_{\g_2} -\frac{\sigma}{2}\1\otimes h_{(-2)}\1$. Then $\omega_{\g, \sigma}$ is a Virasoro element of $V_{k}(\g)$.
\end{proposition}
\pf The proof is similar to that of Proposition \ref{v-B(0,n)}. By Proposition \ref{changeV}, $\omega_{\g, \sigma}$ is a Virasoro element of $V_{k_1}(\g_1)\otimes V_{k_2}(\g_2)$. By Proposition \ref{embedding}, $V_{k}(\g)=V_{k_1}(\g_1)\otimes V_{k_2}(\g_2)\oplus  V^1,$ where $V^1$ is the weak $V_{k_1}(\g_1)\otimes V_{k_2}(\g_2)$-module given in Proposition \ref{embedding}. Moreover, it has been proved in the proof of Theorem \ref{rationalaff} that $V^1$ is an $\omega_{\g, \sigma}$-ordinary module of $V_{k}(\g)=V_{k_1}(\g_1)\otimes V_{k_2}(\g_2)$. This implies that $\omega_{\g, \sigma}$ is a Virasoro element of $V_{k}(\g)$.
\qed

We are now ready to prove the rationality of $V_{k}(\g)$.
\begin{theorem}\label{r-s-2}
Let $\g$ be a simple Lie superalgebra listed in TABLE \ref{tab1}, and $k$ be the corresponding  complex number listed in TABLE \ref{tab1}. Then the vertex superalgebra $V_{k}(\g)$ is $\omega_{\g, \sigma}$-rational.
\end{theorem}
\pf The proof is similar to that of Theorem \ref{r-B(0,n)}. Let $M$ be an $\omega_{\g, \sigma}$-admissible $V_{k}(\g)$-module. By the argument in the proof of Theorem \ref{rationalaff}, $M$ viewed as a weak $V_{k_1}(\g_1)\otimes V_{k_2}(\g_2)$-module belongs to $\mathcal{C}_{\g_1,\g_2}$. Therefore, $M$ belongs to  $\mathcal{D}_{\g,k}$. By Proposition \ref{D-AS}, $M$ is a direct sum of irreducible weak $V_{k}(\g)$-modules. Furthermore, by the similar argument as in the proof of Theorem \ref{r-B(0,n)}, $M$ is a direct sum of irreducible $\omega_{\g, \sigma}$-admissible $V_{k}(\g)$-modules.
\qed

\section{Rationality of certain  minimal $W$-algebras}
In this section, we will show that certain  minimal $W$-algebras are rational with respect to a family of Virasoro elements.
\subsection{Conformal embeddings} In the following, we use $\g=\g_{\bar 0}\oplus \g_{\bar 1}$ to denote a basic classical simple Lie superalgebra. Fix a Cartan subalgebra $\h$  of $\g_{\bar 0}$, and let $\Delta_{\g}$ be the set of roots of $\g$. Fix a minimal root $-\theta$ of $\g$, then we choose root vectors $e_{\theta}$ and $e_{-\theta}$ such that
$$[e_{\theta}, e_{-\theta}]=x\in \h,~~[x, e_{\pm \theta}]=\pm e_{\pm \theta}.$$
Let $(\cdot|\cdot)$ be the  non-degenerate invariant supersymmetric bilinear form of $\g$ such that $(\theta|\theta)=2$. For any complex number $k$, we use $W^k(\g, e_{-\theta})$ to denote the  minimal $W$-algebra of level $k$ defined in \cite{KRW}. It is known that $W^k(\g, e_{-\theta})$  has a unique simple quotient, denoted by $W_k(\g, e_{-\theta})$. In the following, we assume that $k\neq h^{\vee}$, where $h^{\vee}$ denotes the dual Coxeter number of $(\g,(\cdot|\cdot))$. In such a case, it is known that $W^k(\g, e_{-\theta})$  has a Virasoro element $\omega$ defined by the formula (2.2) of \cite{KW3}.

\begin{table}[h]
\centering
\begin{tabular}{|c c c c|}
\hline
$\g$ & $\g^{\natural}$&  $k^{\natural}$& $k$\\
\hline
$G_2$ &  $A_1$& $\frac{1}{2}$ &$-\frac{3}{2}$\\
\hline
$E_6$ &  $A_5$& $-\frac{5}{2}$ &$-\frac{11}{2}$\\
\hline
$E_7$ &  $D_6$& $-\frac{9}{2}$ &$-\frac{17}{2}$\\
\hline
$E_8$ &  $E_7$& $-\frac{17}{2}$ &$-\frac{29}{2}$\\
\hline
$F(4)$ &  $B_3$& $-\frac{13}{4}$ &$\frac{3}{2}$\\
\hline
$C_{n+1}, n\geq 2, (2n+1, 3)=1$ &  $C_n$& $-\frac{4n+5}{6}$& $-\frac{2(n+2)}{3}$\\
\hline
$osp(2|2n), n\geq 3, (2n-1, 3)=1$ &  $D_n$& $\frac{-4n+5}{3}$& $\frac{2(n-2)}{3}$\\
 \hline
 $osp(2|2n+1), n\geq 1, (2n, 3)=1$ &  $B_n$& $\frac{-4n+3}{3}$& $\frac{2n-3}{3}$\\
 \hline
\end{tabular}
\caption{}
\label{tab2}
\end{table}
 For a basic classical simple Lie superalgebra $\g$  listed in TABLE \ref{tab2}, we define the Lie subalgebra $\g^{\natural}$ of $\g$ as in TABLE \ref{tab2}. Let $k, k^{\natural}$ be the complex numbers listed in TABLE \ref{tab2}. It was proved in Theorem 6.8 of \cite{AKMPP2} that the vertex superalgebra $W_k(\g, e_{-\theta})$  has a vertex subalgebra $V_{k^{\natural}}(\g^{\natural})$. Let $\omega_{\g^{\natural}}$  be the canonical Virasoro element of $V_{k^{\natural}}(\g^{\natural})$. It was proved in Theorem 6.8 of \cite{AKMPP2} that $\omega=\omega_{\g^{\natural}}$. For any $v\in W_k(\g, e_{-\theta})$, we denote by $\wt v$ the conformal weight of $v$ under the action of $\omega_{(1)}$. Set $W_k(\g, e_{-\theta})^{\bar 0}=\{v\in W_k(\g, e_{-\theta})|\wt v\in \Z\}$ and $W_k(\g, e_{-\theta})^{\bar 1}=\{v\in W_k(\g, e_{-\theta})|\wt v\in \frac{1}{2}+\Z\}$. Then the following result  has been proved in Theorem 6.8 of \cite{AKMPP2}.
 \begin{theorem}\label{embedding1}
Let $\g$ be a basic classical simple Lie superalgebra listed in TABLE \ref{tab2}, $k, k^{\natural}$ be the complex numbers listed in TABLE \ref{tab2}. Then the subalgebra $W_k(\g, e_{-\theta})^{\bar 0}$ of $W_k(\g, e_{-\theta})$ is isomorphic to $V_{k^{\natural}}(\g^{\natural})$. Moreover, $W_k(\g, e_{-\theta})=W_k(\g, e_{-\theta})^{\bar 0}\oplus W_k(\g, e_{-\theta})^{\bar 1}$, and $W_k(\g, e_{-\theta})^{\bar 1}$ is an irreducible $\omega$-ordinary $V_{k^{\natural}}(\g^{\natural})$-module.
\end{theorem}

\begin{table}[h]
\centering
\begin{tabular}{|c c c c c c|}
\hline
$\g$ & $\g_1$& $\g_2$& $k_1$& $k_2$& $k$\\
\hline
$D_n$ &  $A_1$& $D_{n-2}$& $-1/2$& $\frac{7}{2}-n$& $\frac{3}{2}-n$\\
\hline
$osp(4|2n), n\geq 4, (2n+1, 3)=1$ &  $A_1$& $C_n$& $-1/2$& $-\frac{2n+3}{4}$& $n-\frac{1}{2}$\\
 \hline
\end{tabular}
\caption{}
\label{tab3}
\end{table}
 For a basic classical simple Lie superalgebra $\g$  listed in TABLE \ref{tab3}, we define the Lie subalgebras $\g_1$ and $\g_2$ of $\g$ as in TABLE \ref{tab3}. Let $k, k_1, k_2$ be complex numbers listed in TABLE \ref{tab3}. It was proved in Theorem 6.8 of \cite{AKMPP2} that the vertex superalgebra $W_k(\g, e_{-\theta})$ has a vertex subalgebra $V_{k_1}(\g_1)\otimes V_{k_2}(\g_2)$. Let $\omega_{\g_1}$ and $\omega_{\g_2}$ be the canonical Virasoro elements of $V_{k_1}(\g_1)$ and $V_{k_2}(\g_2)$, respectively. It was proved in Theorem 6.8 of \cite{AKMPP2} that $\omega=\omega_{\g_1}\otimes \1+\1\otimes \omega_{\g_2}$. Set $W_k(\g, e_{-\theta})^{\bar 0}=\{v\in W_k(\g, e_{-\theta})|\wt v\in \Z\}$ and $W_k(\g, e_{-\theta})^{\bar 1}=\{v\in W_k(\g, e_{-\theta})|\wt v\in \frac{1}{2}+\Z\}$. Then the following result  has been proved in Theorem 6.8 of \cite{AKMPP2}.
\begin{theorem}\label{embedding2}
Let $\g$ be a basic classical simple Lie superalgebra listed in TABLE \ref{tab3}, $k, k_1, k_2$ be the complex numbers listed in TABLE \ref{tab3}. Then the subalgebra $W_k(\g, e_{-\theta})^{\bar 0}$ of $W_k(\g, e_{-\theta})$ is isomorphic to  $V_{k_1}(\g_1)\otimes V_{k_2}(\g_2)$. Moreover, $W_k(\g, e_{-\theta})=W_k(\g, e_{-\theta})^{\bar 0}\oplus W_k(\g, e_{-\theta})^{\bar 1}$, and $W_k(\g, e_{-\theta})^{\bar 1}$ is an irreducible $\omega$-ordinary $V_{k_1}(\g_1)\otimes V_{k_2}(\g_2)$-module.
\end{theorem}
\subsection{Rationality of  minimal $W$-algebras at conformal level} In this subsection, we show that the  minimal $W$-algebras in Theorems \ref{embedding1}, \ref{embedding2} are rational with respect to a family of Virasoro elements. We will prove the rationality by using Theorem \ref{scextension}. First, let $\g$ be a basic classical simple Lie superalgebra listed in TABLE \ref{tab2}, $\g^{\natural}$ be the Lie subalgebra of $\g$ listed in TABLE \ref{tab2}, $k, k^{\natural}$ be complex numbers listed in TABLE \ref{tab2}. Then the subalgebra $W_k(\g, e_{-\theta})^{\bar 0}$ of $W_k(\g, e_{-\theta})$ is isomorphic to $V_{k^{\natural}}(\g^{\natural})$. Let $\mathcal{C}_{\g^{\natural},k^{\natural}}$ be the subcategory of weak $V_{k^{\natural}}(\g^{\natural})$-module category defined in subsection \ref{admissiblesec}. We define $\mathcal{E}_{\g,k}$ to be the subcategory of weak $W_k(\g, e_{-\theta})$-module category such that $M$ is an object in $\mathcal{E}_{\g,k}$ if and only if $M$ viewed as a weak $V_{k^{\natural}}(\g^{\natural})$-module belongs to  $\mathcal{C}_{\g^{\natural},k^{\natural}}$. Then we have the following result.
\begin{proposition}\label{E-W-1}
Let $\mathcal{E}_{\g,k}$ be the category defined above. Then $\mathcal{E}_{\g,k}$ is semisimple.
\end{proposition}
\pf
Let $\g$ be a basic classical simple Lie superalgebra listed in TABLE \ref{tab2}, $k, k^{\natural}$ be the complex numbers listed in TABLE \ref{tab2}.  Then $k^{\natural}$ is an admissible level of $\hat{\g^{\natural}}$. By Theorem \ref{semisimple},  $\mathcal{C}_{\g^{\natural},k^{\natural}}$ is semisimple.

By Theorem \ref{embedding1}, the subalgebra $W_k(\g, e_{-\theta})^{\bar 0}$ of $W_k(\g, e_{-\theta})$ is isomorphic to  $V_{k^{\natural}}(\g^{\natural})$. Moreover, $W_k(\g, e_{-\theta})=W_k(\g, e_{-\theta})^{\bar 0}\oplus W_k(\g, e_{-\theta})^{\bar 1}$ and $W_k(\g, e_{-\theta})^{\bar 1}$ is an irreducible $\omega$-ordinary $V_{k^{\natural}}(\g^{\natural})$-module. Set $$W_k(\g, e_{-\theta})^{\bar 1}\cdot W_k(\g, e_{-\theta})^{\bar 1}=\{u_{(n)}v|u, v\in W_k(\g, e_{-\theta})^{\bar 1}, n\in \Z\}.$$ By the definitions of $W_k(\g, e_{-\theta})^{\bar 0}$ and $W_k(\g, e_{-\theta})^{\bar 1}$, we have $W_k(\g, e_{-\theta})^{\bar 1}\cdot W_k(\g, e_{-\theta})^{\bar 1}=W_k(\g, e_{-\theta})^{\bar 0}$. It follows from Theorem \ref{scextension} that $\mathcal{E}_{\g,k}$ is semisimple.
\qed

We now show that  $W_k(\g, e_{-\theta})$ has a family of Virasoro elements.
\begin{proposition}\label{V-W-1}
Let $\g, \g^{\natural}, k, k^{\natural}$ be as above, $h$ be an element in $C^{\circ}_{\g^{\natural}, -}$, $\sigma$ be a rational number such that $-1<\frac{\sigma}{2}\alpha(h)<0$ for all $\alpha\in \Delta_{\g^{\natural}}^+$, $\omega_{\g^{\natural}}$ be the canonical Virasoro element of $V_{k^{\natural}}(\g^{\natural})$ defined as in Theorem \ref{affine}. Set $\omega_{\g^{\natural}, \sigma}=\omega_{\g^{\natural}}-\frac{\sigma}{2}h_{(-2)}\1$. Then $\omega_{\g^{\natural}, \sigma}$ is a Virasoro element of $W_k(\g, e_{-\theta})$.
\end{proposition}
\pf By Proposition \ref{changeV}, $\omega_{\g^{\natural}, \sigma}$ is a Virasoro element of $V_{k^{\natural}}(\g^{\natural})$. By Theorem \ref{embedding1},  the subalgebra $W_k(\g, e_{-\theta})^{\bar 0}$ of $W_k(\g, e_{-\theta})$ is isomorphic to  $V_{k^{\natural}}(\g^{\natural})$. Moreover, $W_k(\g, e_{-\theta})=W_k(\g, e_{-\theta})^{\bar 0}\oplus W_k(\g, e_{-\theta})^{\bar 1}$ and $W_k(\g, e_{-\theta})^{\bar 1}$ is an irreducible $\omega$-ordinary $V_{k^{\natural}}(\g^{\natural})$-module. The explicit forms of $W_k(\g, e_{-\theta})^{\bar 1}$ as an $\omega$-ordinary $V_{k^{\natural}}(\g^{\natural})$-module were given in Theorem 6.8 of \cite{AKMPP2}. In particular, it was proved in  Theorem 6.8 of \cite{AKMPP2} that $W_k(\g, e_{-\theta})^{\bar 1}$ is an irreducible highest weight module of $\hat{\g^{\natural}}$. Thus, it follows from the proof of Theorem \ref{rationalaff} that $W_k(\g, e_{-\theta})^{\bar 1}$ is an $\omega_{\g^{\natural}, \sigma}$-ordinary module of $V_{k^{\natural}}(\g^{\natural})$. This implies that $\omega_{\g^{\natural}, \sigma}$ is a Virasoro element of $W_k(\g, e_{-\theta})$.
\qed

We now prove that the vertex superalgebras $W_k(\g, e_{-\theta})$ in Theorem \ref{embedding1} are $\omega_{\g^{\natural}, \sigma}$-rational.
\begin{theorem}\label{r-w-1}
Let $\g$ be a basic classical simple Lie superalgebra listed in TABLE \ref{tab2}, $k$ be the complex number listed in TABLE \ref{tab2}. Then the vertex superalgebra $W_k(\g, e_{-\theta})$  is $\omega_{\g^{\natural}, \sigma}$-rational.
\end{theorem}
\pf Let $M$ be an $\omega_{\g^{\natural}, \sigma}$-admissible $W_k(\g, e_{-\theta})$-module.  By the argument in the proof of Theorem \ref{rationalaff}, $M$ viewed as a weak $V_{k^{\natural}}(\g^{\natural})$-module belongs to $\mathcal{C}_{\g^{\natural},k^{\natural}}$. Therefore, $M$ belongs to  $\mathcal{E}_{\g,k}$. By Proposition \ref{E-W-1}, $M$ is a direct sum of irreducible weak $W_k(\g, e_{-\theta})$-modules. Furthermore, by the similar argument as in the proof of Theorem \ref{r-B(0,n)}, $M$ is a direct sum of irreducible $\omega_{\g^{\natural}, \sigma}$-admissible $W_k(\g, e_{-\theta})$-modules.
\qed

\vskip.25cm
Let $\g$ be a basic classical simple Lie superalgebra listed in TABLE \ref{tab3}, $\g_1$, $\g_2$ be the Lie subalgebras of $\g$ listed in TABLE \ref{tab3}, $k, k_1, k_2$ be the complex numbers listed in TABLE \ref{tab3}. Then the subalgebra $W_k(\g, e_{-\theta})^{\bar 0}$ of $W_k(\g, e_{-\theta})$ is isomorphic to  $V_{k_1}(\g_1)\otimes V_{k_2}(\g_2)$. Let $\mathcal{C}_{\g_1, \g_2}$ be the subcategory of weak $V_{k_1}(\g_1)\otimes V_{k_2}(\g_2)$-module category such that $M$ is an object in $\mathcal{C}_{\g_1, \g_2}$ if and only if $M$ satisfies the following two conditions:\\
(1) $M$ viewed as a weak $V_{k_1}(\g_1)$-module belongs to $\mathcal{C}_{\g_1, k_1}$.\\
(2) $M$ viewed as a weak $V_{k_2}(\g_2)$-module belongs to $\mathcal{C}_{\g_2, k_2}$.\\
Then we have the following result.
\begin{proposition}\label{C-W-2}
Let $\mathcal{C}_{\g_1, \g_2}$ be the category defined above. Then $\mathcal{C}_{\g_1, \g_2}$ is semisimple.
\end{proposition}
\pf Note that $k_1, k_2$ are admissible levels of $\hat{\g}_1, \hat{\g}_2$, respectively. By the similar argument as in the proof of Proposition \ref{semisimpleten1}, $\mathcal{C}_{\g_1, \g_2}$ is semisimple.
\qed

\vskip.25cm
Let $\g$ be a basic classical simple Lie superalgebra listed in TABLE \ref{tab3}, $\g_1$, $\g_2$ be the Lie subalgebras of $\g$ listed in TABLE \ref{tab3}, $k, k_1, k_2$ be the complex numbers listed in TABLE \ref{tab3}. We define $\mathcal{E}_{\g,k}$ to be the subcategory of weak $W_k(\g, e_{-\theta})$-module category such that $M$ is an object in $\mathcal{E}_{\g,k}$ if and only if $M$ viewed as a weak $V_{k_1}(\g_1)\otimes V_{k_2}(\g_2)$-module belongs to  $\mathcal{C}_{\g_1, \g_2}$. Then we have the following result.
\begin{proposition}\label{E-W-2}
Let $\mathcal{E}_{\g,k}$ be the category defined as above. Then $\mathcal{E}_{\g,k}$ is semisimple.
\end{proposition}
\pf The proof is similar to that of Proposition \ref{E-W-1}.
Let $\g$ be a basic classical simple Lie superalgebra listed in TABLE \ref{tab3}, $\g_1$, $\g_2$ be the Lie subalgebras of $\g$ listed in TABLE \ref{tab3}, $k, k_1, k_2$ be the complex numbers listed in TABLE \ref{tab3}.
By Theorem \ref{embedding2}, the subalgebra $W_k(\g, e_{-\theta})^{\bar 0}$ of $W_k(\g, e_{-\theta})$ is isomorphic to  $V_{k_1}(\g_1)\otimes V_{k_2}(\g_2)$. Moreover, $W_k(\g, e_{-\theta})=W_k(\g, e_{-\theta})^{\bar 0}\oplus W_k(\g, e_{-\theta})^{\bar 1}$ and $W_k(\g, e_{-\theta})^{\bar 1}$ is an irreducible $\omega$-ordinary $V_{k_1}(\g_1)\otimes V_{k_2}(\g_2)$-module. Set $$W_k(\g, e_{-\theta})^{\bar 1}\cdot W_k(\g, e_{-\theta})^{\bar 1}=\{u_{(n)}v|u, v\in W_k(\g, e_{-\theta})^{\bar 1}, n\in \Z\}.$$ By the definitions of $W_k(\g, e_{-\theta})^{\bar 0}$ and $W_k(\g, e_{-\theta})^{\bar 1}$, we have $W_k(\g, e_{-\theta})^{\bar 1}\cdot W_k(\g, e_{-\theta})^{\bar 1}=W_k(\g, e_{-\theta})^{\bar 0}$. It follows from Proposition \ref{C-W-2} and Theorem \ref{scextension} that $\mathcal{E}_{\g,k}$ is semisimple.
\qed

We now show that  $W_k(\g, e_{-\theta})$ has a family of Virasoro elements.
\begin{proposition}\label{V-W-2}
Let $\g$, $\g_1$, $\g_2$, $k, k_1, k_2$ be as above. For $i=1, 2$, let   $h_i $ be an element in $C^{\circ}_{\g_i, -}$, $\sigma_i$ be a rational number such that $-1<\frac{\sigma_i}{2}\alpha(h)<0$ for all $\alpha\in \Delta_{\g_i}^+$, $\omega_{\g_i}$ be the canonical Virasoro element of $V_{k_i}(\g_i)$ defined as in Theorem \ref{affine}. Set $\omega_{\g, \sigma_1, \sigma_2}=\omega_{\g_1}\otimes \1+1\otimes \omega_{\g_2} -\frac{\sigma_1}{2}(h_1)_{(-2)}\1\otimes\1 -\frac{\sigma_2}{2}\1\otimes (h_2)_{(-2)}\1$. Then $\omega_{\g, \sigma_1, \sigma_2}$ is a Virasoro element of $W_k(\g, e_{-\theta})$.
\end{proposition}
\pf We first consider the case that $\g=D_n$. Then $\g_1=A_1$, $\g_2=D_{n-2}$. By Proposition \ref{changeV}, $\omega_{A_1, \sigma_1}$ is a Virasoro element of $V_{-1/2}(A_1)$, and $\omega_{D_{n-2}, \sigma_2}$ is a Virasoro element of $V_{\frac{7}{2}-n}(D_{n-2})$. By Theorem \ref{embedding1},  the subalgebra $W_{\frac{3}{2}-n}(D_n, e_{-\theta})^{\bar 0}$ of $W_{\frac{3}{2}-n}(D_n, e_{-\theta})$ is isomorphic to  $V_{-1/2}(A_1)\otimes V_{\frac{7}{2}-n}(D_{n-2})$. Moreover, $$W_{\frac{3}{2}-n}(D_n, e_{-\theta})=W_{\frac{3}{2}-n}(D_n, e_{-\theta})^{\bar 0}\oplus W_{\frac{3}{2}-n}(D_n, e_{-\theta})^{\bar 1}$$ and $W_{\frac{3}{2}-n}(D_n, e_{-\theta})^{\bar 1}$ is an irreducible $\omega$-ordinary $V_{-1/2}(A_1)\otimes V_{\frac{7}{2}-n}(D_{n-2})$-module isomorphic to $L_{A_1}(-\frac{1}{2},\Lambda_1)\otimes L_{D_{n-2}}(\frac{7}{2}-n,\Lambda_1)$.  By the proof of Theorem \ref{rationalaff}, $L_{A_1}(-\frac{1}{2},\Lambda_1)\otimes L_{D_{n-2}}(\frac{7}{2}-n,\Lambda_1)$ is an $\omega_{D_n, \sigma_1, \sigma_2}$-ordinary module of $V_{-1/2}(A_1)\otimes V_{\frac{7}{2}-n}(D_{n-2})$. This implies that $\omega_{D_n, \sigma_1, \sigma_2}$ is a Virasoro element of $W_{\frac{3}{2}-n}(D_n, e_{-\theta})$. Similarly, $\omega_{osp(4|2n), \sigma_1, \sigma_2}$ is a Virasoro element of $W_{n-\frac{1}{2}}(osp(4|2n), e_{-\theta})$.
\qed

We now prove that the vertex superalgebras $W_k(\g, e_{-\theta})$ in Theorem \ref{embedding2} are $\omega_{\g, \sigma_1, \sigma_2}$-rational.
\begin{theorem}\label{r-w-2}
Let $\g$ be a basic classical simple Lie superalgebra listed in TABLE \ref{tab3}, $k$ be the complex number listed in TABLE \ref{tab3}. Then the vertex superalgebra $W_k(\g, e_{-\theta})$  is $\omega_{\g, \sigma_1, \sigma_2}$-rational.
\end{theorem}
\pf Let $\g$ be a basic classical simple Lie superalgebra listed in TABLE \ref{tab3}, $k$ be the complex number listed in TABLE \ref{tab3}, $M$ be an $\omega_{\g, \sigma_1, \sigma_2}$-admissible $W_k(\g, e_{-\theta})$-module.  By the argument in the proof of Theorem \ref{rationalaff}, $M$ viewed as a weak $V_{k_1}(\g_1)\otimes V_{k_2}(\g_2)$-module belongs to $\mathcal{C}_{\g_1,\g_2}$. Therefore, $M$ belongs to  $\mathcal{E}_{\g,k}$. By Proposition \ref{E-W-2}, $M$ is a direct sum of irreducible weak $W_k(\g, e_{-\theta})$-modules. Furthermore, by the similar argument as in the proof of Theorem \ref{r-B(0,n)}, $M$ is a direct sum of irreducible $\omega_{\g, \sigma_1, \sigma_2}$-admissible $W_k(\g, e_{-\theta})$-modules.
\qed

\vskip.25cm
{\bf Acknowledgement}. The author wish to thank Professor Tomoyuki Arakawa for  suggestions.

\end{document}